\documentclass[12pt,reqno]{amsart}
\usepackage{amsmath}
\usepackage{amssymb}
\usepackage{verbatim}
\usepackage{mathrsfs}
\usepackage{graphicx} 
\usepackage{caption}
\usepackage{subcaption}
\usepackage{epstopdf}
\allowdisplaybreaks

\usepackage[toc]{appendix}
\usepackage{wavelet}

\usepackage[top=0.65in,bottom=0.64in,left=0.7in,right=0.7in]{geometry}

\newtheorem{lemma}{Lemma}

\newtheorem{theorem}[lemma]{Theorem}

\newtheorem{example}{Example}

\newtheorem{algorithm}{Algorithm}

\numberwithin{equation}{section}
\numberwithin{lemma}{section}

\begin{document}

\title[Wavelet Methods for Differential Equations (ICCM 2019)]{Wavelet-based Methods for Numerical Solutions of Differential Equations}

\author{Bin Han, Michelle Michelle, and Yau Shu Wong}

\thanks{Research supported in part by
Natural Sciences and Engineering Research Council (NSERC) of Canada and Alberta Innovates}

\address{Department of Mathematical and Statistical Sciences,
University of Alberta, Edmonton,\quad Alberta, Canada T6G 2G1.
\quad {\tt bhan@ualberta.ca}\quad {\tt mmichell@ualberta.ca} \quad {\tt yauwong@ualberta.ca}
}

\makeatletter \@addtoreset{equation}{section} \makeatother
\begin{abstract}
Wavelet theory has been well studied in recent decades. Due to their appealing features such as sparse multiscale representation and fast algorithms, wavelets have enjoyed many tremendous successes in the areas of signal/image processing and computational mathematics. This paper primarily intends to shed some light on the advantages of using wavelets in the context of numerical differential equations. We shall identify a few prominent problems in this field and recapitulate some important results along these directions.
Wavelet-based methods for numerical differential equations offer the advantages of sparse matrices with uniformly bounded small condition numbers.
We shall demonstrate wavelets' ability in solving some one-dimensional differential equations: the biharmonic equation and the Helmholtz equation with high wave numbers (of magnitude $\bo(10^{4})$ or larger).
\end{abstract}


\maketitle

\pagenumbering{arabic}

\section{Introduction and Motivations}
Being unconditional bases in many function spaces,
wavelets are sparse multiscale representation systems, which serve as excellent approximation tools for various types of functions and signals. The theory of wavelet analysis can be traced back as early as Haar's discovery in 1909, which is now known as the Haar orthogonal wavelet system. It is the simplest example of a wavelet, since it takes the form of an indicator function supported on the unit interval. However, wavelet analysis only experienced a surge of interest and rapid development after the emergence of
Morlet \cite{ggm84}, Meyer \cite{m90}, and Daubechies \cite{d88} wavelets as well as the notion of multiresolution analysis \cite{m89,m90}. At the present time, the theory of wavelet analysis is considerably well-established (e.g., see \cite{hanbook}) and its applications are far-reaching. In the context of computational mathematics, wavelets have demonstrated their ability in tackling problems related to signal/image processing and numerical differential equations (DEs). The main focus of this paper is to explore the upper hand we gain from utilizing wavelets in the context of numerical DEs. We refer interested readers to the recent book \cite{hanbook} and references therein for a similar discussion in the context of signal and image processing.

\begin{figure}[htbp]
	\centering
	\begin{subfigure}[b]{0.28\textwidth}
		 \includegraphics[width=\textwidth]{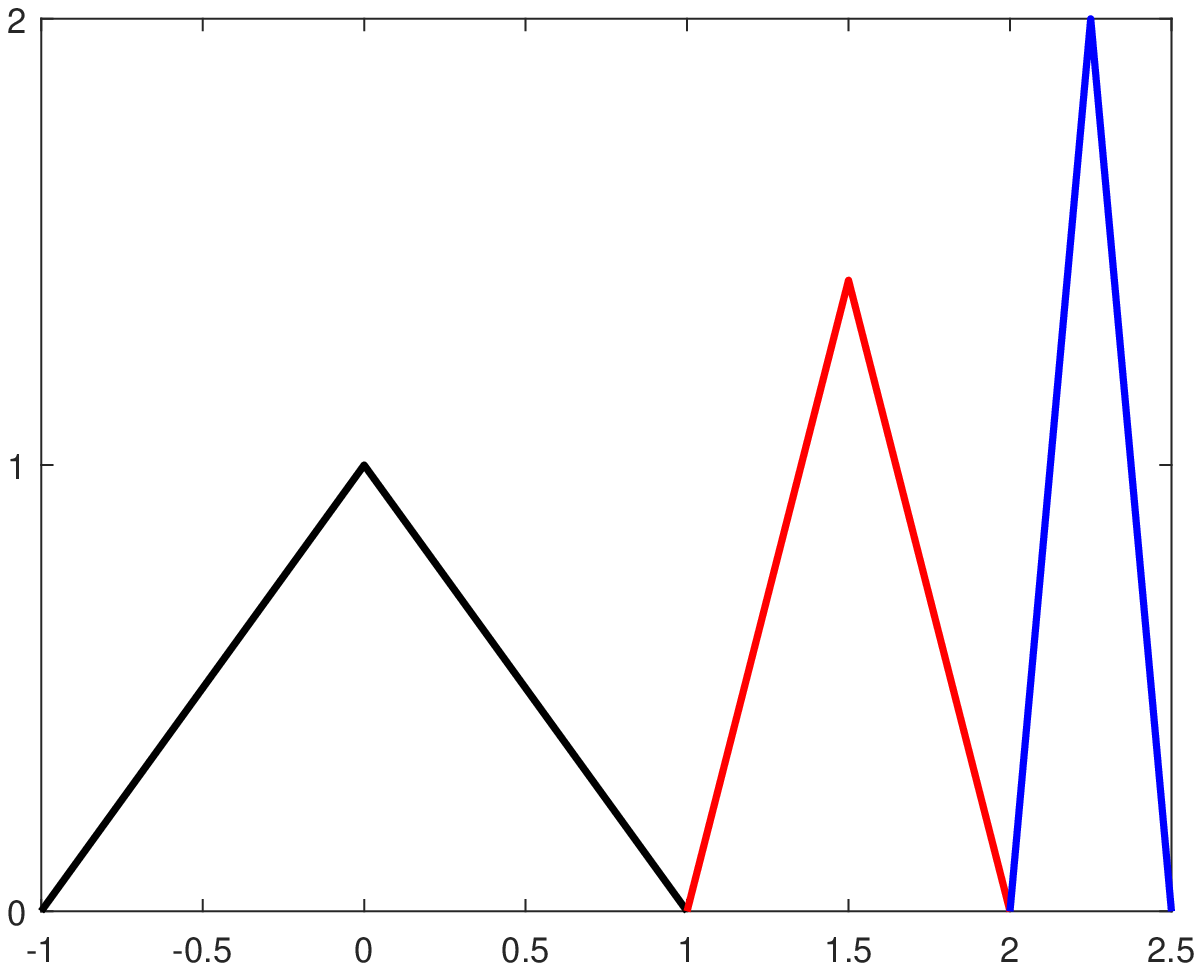}
	\end{subfigure}
	\begin{subfigure}[b]{0.28\textwidth}
		 \includegraphics[width=\textwidth]{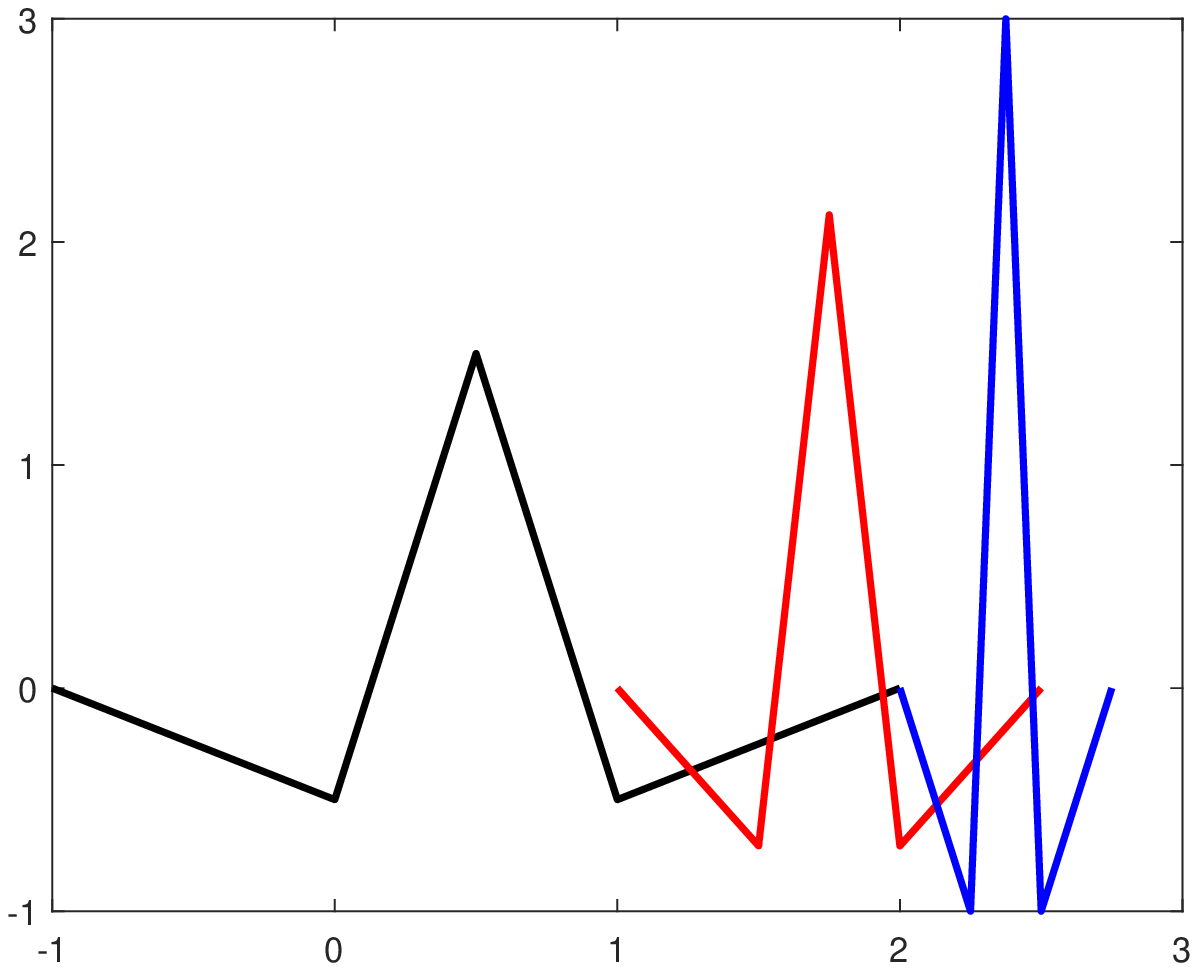}
	\end{subfigure}
	\caption{Sample elements from a Riesz wavelet basis. Left: $\phi$ (black), $2^{1/2}\phi(2\cdot-3)$ (red), $2\phi(2^2\cdot-9)$ (blue). Right: $\psi$ (black), $2^{1/2}\psi(2\cdot-3)$ (red), $2\psi(2^2\cdot-9)$ (blue).}
	\label{fig:phipsishifts}
\end{figure}

To elucidate the concept of wavelets, we recall some basic facts and present some examples. Define $\NN:=\N\cup\{0\}$ and let $m \in \NN$. Recall that the Sobolev space $H^{m}(\R)$ contains all functions $f$ on $\R$ such that $f, f', \dots, f^{(m-1)}$ are absolutely continuous on $\R$ and $f, \dots, f^{(m)} \in L_{2}(\R)$. When $m=0$, we naturally have $H^{0}(\R) = L_{2}(\R)$. Let $\phi=(\phi_{1},\dots,\phi_{r})^{\tp}$ and $\psi=(\psi_{1},\dots,\psi_{s})^{\tp}$ be vectors of tempered distributions or functions on $\R$. For $J\in \Z$, a \emph{wavelet affine system} for $H^{m}(\R)$ is defined as
\begin{equation} \label{ASR}
\begin{split}
\AS_J^{m}(\phi;\psi):=
&\{2^{J(1/2-m)}\phi_\ell(2^J \cdot-k) \setsp k\in \Z, 1\le \ell\le r\} \\
&\quad \cup
\{2^{j(1/2-m)}\psi_\ell(2^j\cdot-k) \setsp j\ge J, k\in \Z, 1\le \ell\le s\}.
\end{split}
\end{equation}
We define $\AS_J(\phi;\psi):=\AS_J^0(\phi;\psi)$ for $m=0$. We often refer $j$ and $k$ in \eqref{ASR} as the scale/resolution level and the integer shift, respectively. Hence, a wavelet affine system is simply a set that contains dilated and translated versions of functions in all the entries of $\phi$ and $\psi$. We shall see later that the wavelet function $\psi$ is also generated from a linear combination of dilated and translated versions of the refinable (vector) function $\phi$. Note that $\phi$ in Figure~\ref{fig:phipsishifts} is the hat function or B-spline of order 2. In fact, B-splines and Hermite splines are two popular examples of refinable functions, $\phi$, due to their analytic expressions and wide usage in computational mathematics.
The wavelet affine system $\AS_J^{m}(\phi;\psi)$ is a Riesz basis for $\HH{m}$ if the following two conditions are satisfied: (1) the linear span of $\AS_J^{m}(\phi;\psi)$ is dense in $\HH{m}$, and (2) there exist positive constants $C_1$ and $C_2$ such that
{\small
\begin{align*}
C_1 \Big(\sum_{\ell=1}^r \sum_{k\in \Z} |v_{\ell,k}|^2
&+\sum_{j=J}^\infty \sum_{\ell=1}^s \sum_{k\in \Z}
|w_{\ell,j;k}|^2\Big) \le
\Big\| \sum_{\ell=1}^r \sum_{k\in \Z} v_{\ell,k} 2^{J(1/2-m)}\phi_\ell(2^J \cdot-k)\\
&+ \sum_{j=J}^\infty \sum_{\ell=1}^s \sum_{k\in \Z} w_{\ell,j;k} 2^{j(1/2-m)} \psi_\ell(2^j\cdot-k)\Big\|^2_{\HH{m}}
\le C_2 \Big( \sum_{\ell=1}^r \sum_{k\in \Z} |v_{\ell,k}|^2
+ \sum_{j=J}^\infty \sum_{\ell=1}^s \sum_{k\in \Z} |w_{\ell,j;k}|^2\Big)
\end{align*}
}
for all finitely supported sequences $\{ v_{\ell,k}\}_{1\le \ell\le r,k\in \Z}$ and $\{w_{\ell,j;k}\}_{1\le \ell\le s,j\ge J, k\in \Z}$. The best possible constants $C_{1}$ and $C_{2}$ are respectively called the \emph{lower and upper Riesz bounds} of $\AS_J^{m}(\phi;\psi)$.
It is known in \cite{han10,hanbook} that $\AS_J^{m}(\phi;\psi)$ is a Riesz basis for $H^m(\R)$ if and only if it is a Riesz basis for $H^m(\R)$ for all $J\in \Z$.
Hence, we call $\{\phi;\psi\}$ \emph{a Riesz wavelet in the Sobolev space} $\HH{m}$ if $\AS_0^{m}(\phi;\psi)$ is a Riesz basis for $H^m(\R)$.
Moreover, the ratio $C_2/C_1$ of its Riesz bounds is called the \emph{condition number} of the Riesz wavelet $\{\phi;\psi\}$ (or of the Riesz wavelet basis
$\AS_0^{m}(\phi;\psi)$) in the Sobolev space $\HH{m}$.
See Figure~\ref{fig:phipsishifts} for sample elements taken from such Riesz wavelet bases.
If $r=1$, we often call $\{\phi;\psi\}$ a scalar Riesz wavelet. Meanwhile, if $r>1$, we often call $\{\phi;\psi\}$ a Riesz multiwavelet due to the fact that $\phi$ is a vector function. For the sake of convenience, we simply use wavelets or Riesz wavelets to refer to both scalar Riesz wavelets and Riesz multiwavelets.

We now discuss why it is favourable to employ wavelets in solving numerical DEs. Wavelets are commonly used in the finite element/Galerkin method. Firstly, we have a well-established theory to guarantee that they completely span the space of interest. Thus, we can confidently conclude that the approximated solution approaches the true solution as the relative error decreases. Secondly, Riesz wavelets are well-conditioned bases, which give rise to coefficient matrices with uniformly bounded (possibly) small condition numbers. This property is essential for fast convergence in numerical schemes. Moreover, one standard technique to handle multidimensional DEs is to take the tensor product of a univariate Riesz wavelet with a small condition number. This will prevent the condition number from growing too rapidly with respect to the dimension. In the finite element method, we typically use the hat function as our basis. We shall demonstrate (see Table~\ref{table:cnlegall}) that the condition number of a wavelet basis is far smaller than that of the hat functions (i.e., the standard finite element method), even though both of them generate the same space. In fact, the later condition number seems to go unbounded as the scale level increases. Thirdly, there is a lot of flexibility in designing wavelets. The speed of convergence in the continuous Galerkin method hinges on the polynomial reproduction order of the employed basis. Hence, we can control the speed by simply picking a (primal) refinable function with a suitable polynomial reproduction order. Given this choice, we can construct a (dual) refinable function having a prescribed polynomial reproduction order, which in turn affects the wavelet's vanishing moments. The latter plays an indispensable role in enhancing the sparsity of the coefficient matrix. Typically, we may even end up with some freedom to construct our wavelet basis so that ratio of the Riesz bounds (i.e., the condition number of the Riesz wavelet basis) is as small as possible. Note that high polynomial reproduction order often comes at the cost of longer supports. However, we have the ability to strike the right balance depending on the problem in hand. Once a wavelet basis is aptly chosen, we can use a suitable construction procedure (e.g., \cite{hmw19}) to ensure the boundary conditions (e.g., Dirichlet, Neumann, etc.) are satisfied. In some instances, wavelets can be constructed so that in order to obtain the numerical solution, no linear system needs to be solved. The fourth advantage that wavelets offer comes in the assembly of coefficient matrices in the continuous Galerkin method. In some cases, there is effectively no need to use quadrature to calculate the inner products of the basis functions. Such information can be efficiently obtained from the eigenvector of the transition operator and exploiting the refinability structure of our wavelet basis (i.e., via the fast wavelet transform). Furthermore, we can design an efficient quadrature to compute the inner products of our wavelet basis and the source term as in \cite[Lemma 7.5.6]{hanbook}.

\begin{figure}[htbp]
	\centering
	\begin{subfigure}[b]{0.28\textwidth}
		 \includegraphics[width=\textwidth]{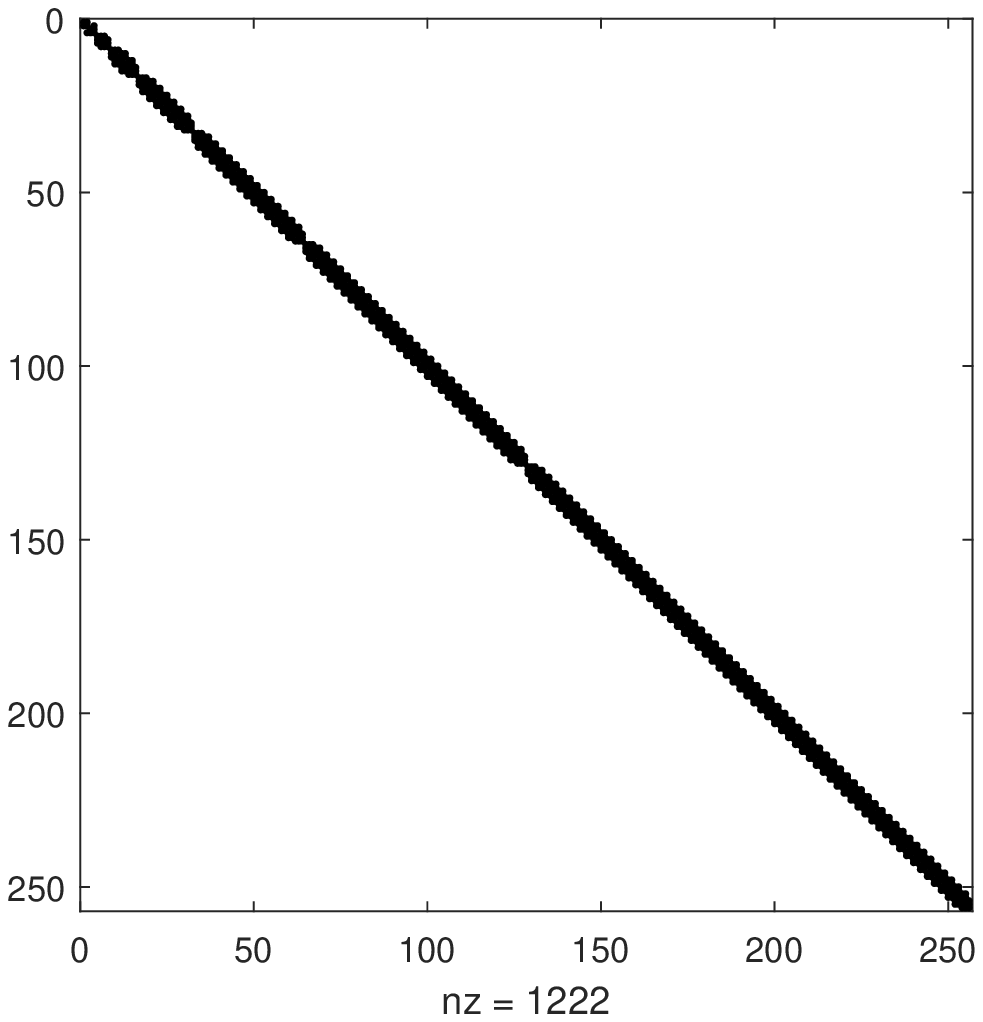}
	\end{subfigure}
	\begin{subfigure}[b]{0.28\textwidth}
		 \includegraphics[width=\textwidth]{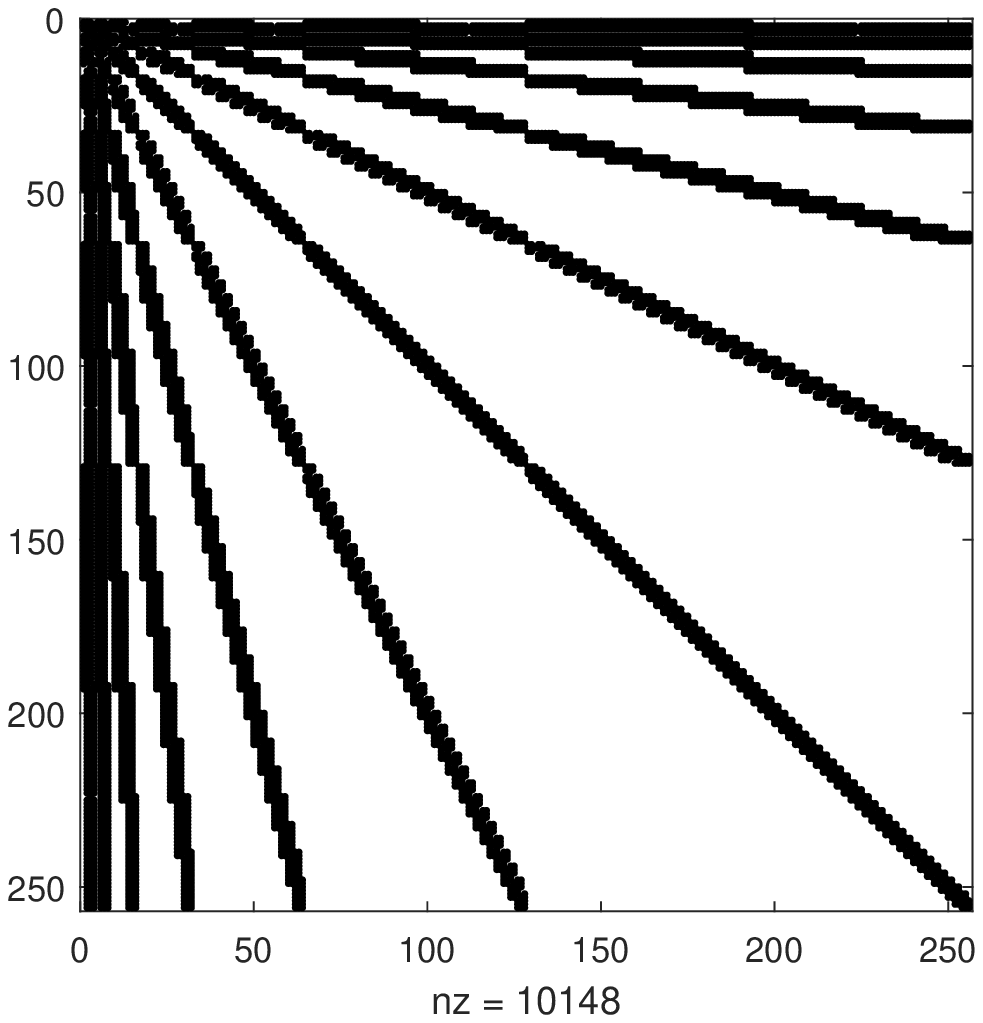}
	\end{subfigure}
	\caption{Left: the sparsity structure of a stiffness matrix stemming from an $m$-th order derivative orthogonal wavelet. Right: the sparsity structure of a mass matrix stemming from an $m$-th order derivative orthogonal wavelet.}
	\label{fig:coefmat}
\end{figure}

There are two fundamental problems in wavelet-based methods for numerical DEs, which act as the focal points of our discussion. The first one is in the construction of suitable wavelets to handle the problem in hand. Having a sparse coefficient matrix is very much desired in solving numerical DEs. Motivated by the model problem
\begin{equation} \label{dow:de}
u^{(2m)}(x)+ \alpha u(x)=f(x), \quad x \in \mathcal{I},
\end{equation}
where $\alpha \in \R$ and $\mathcal{I}$ is a bounded interval on $\R$, \cite{hm19} provides the necessary and sufficient conditions for the construction of wavelets whose $m$-th order derivatives are orthogonal. Prior to \cite{hm19}, these conditions had been unknown. However, examples of such wavelets have long been known in the literature \cite{cr00, jl06, jz11}. Recall that a Riesz wavelet $\{\phi;\psi\}$ in the Sobolev space $H^{m}(\R)$ is called \emph{an $m$-th order derivative orthogonal Riesz wavelet in} $H^{m}(\R)$ if it satisfies the following two properties:
\begin{equation} \label{modo:cond:phi}
\la \psi^{(m)}, \phi^{(m)}(\cdot-k)\ra=0,\qquad \forall\, k\in \Z,
\end{equation}
and
\begin{equation} \label{modo:cond:psi}
\la \psi^{(m)}(2^j\cdot-k), \psi^{(m)}(2^{j'}\cdot-k')\ra=0,\qquad \forall\, k,k'\in \Z, j,j'\in \NN \quad \mbox{with}\; j\ne j'.
\end{equation}
It is not hard to see that the weak formulation of \eqref{dow:de} is
\begin{equation*}
\langle u^{(2m)},v \rangle
+\alpha\la u,v\ra
= \sum_{k=0}^{m-1} (-1)^{k} \left.\left(u^{(2m-k-1)} \ol{v^{(k)}}\right)\right\lvert_{\mathcal{I}} + (-1)^{m}\langle u^{(m)},v^{(m)}\rangle + \alpha \langle u,v \rangle = \langle f,v \rangle
\end{equation*}
for all $v \in H^{m}(\mathcal{I})$. Indeed, employing an $m$-th order derivative-orthogonal Riesz wavelets gives rise to a stiffness matrix with a nice block diagonal matrix structure. On the other hand, the corresponding mass matrix has a sparse finger-like structure. See Figure \ref{fig:coefmat}. Let $\kappa$ be the condition number of a matrix; i.e., $\kappa$ is the ratio of the largest and smallest singular values. For a symmetric positive definite matrix, $\kappa$ is simply the ratio of the largest and smallest eigenvalues. The contribution of $\kappa((-1)^{m}\langle u^{(m)},v^{(m)}\rangle)$ in the condition number of the coefficient matrix, $\kappa((-1)^{m}\langle u^{(m)},v^{(m)}\rangle + \alpha \langle u,v \rangle)$, often overpowers that of $\kappa(\langle u,v \rangle)$. Hence, $\kappa((-1)^{m}\langle u^{(m)},v^{(m)}\rangle + \alpha \langle u,v \rangle)$ is typically small (i.e., well-conditioned). We shall revisit this result in Section 2.
In the same vein, \cite{cf16, ds10} impose an orthogonality condition on the first derivative of the wavelets (generated from Hermite cubic splines) to achieve better sparsity. However, their methods differ from \cite{hm19} in that they substitute the wavelet functions at even integers with a pair that have a longer support.

Just like in many other applications, numerical DEs are typically defined on a bounded domain. For rectangular bounded domains, one strategy is to take the tensor product of wavelets constructed on a bounded interval. Hence, one critical question that we need to resolve is how to optimally construct wavelets on a bounded interval. This serves as the second fundamental problem. Our starting point is a Riesz wavelet on the real line. Ideally, we want to have a tractable algorithm at our disposal that gives us boundary elements with simple structure and allows us to retain all desirable properties that our starting wavelet has. Additionally, we want to ensure that the wavelet formed on the interval of interest indeed spans the whole space. This is in fact a long-standing problem in wavelet analysis. The fact that the analysis of multiwavelets is far more intricate than their scalar counterparts just adds to the difficulty. Nevertheless, any construction we propose needs to be applicable to multiwavelets, since they have a key advantage over scalar wavelets in that they generally have higher smoothness and vanishing moments for a given support. In other words, multiwavelets may attain the desired smoothness and vanishing moments with a much shorter support. Having a short support is indeed necessary to produce boundary elements with simple structure. Tremendous research efforts have gone into developing various construction procedures, but virtually all of them suffer from some serious shortcomings in one form or another. It is not until recently that \cite{hmw19} presents a general construction that gives rise to all possible square integrable compactly supported biorthogonal wavelets on a semi-infinite or bounded interval satisfying maximum vanishing moments and/or boundary conditions. Many existing constructions in the literature serve as special cases of the construction proposed in \cite{hmw19}. We shall elaborate on this further in Section 3.

Next, we mention other equally significant studies in the topic of wavelet-based methods for numerical DEs. First off, the book \cite{cohenbook} provides a comprehensive treatment of the subject. Some key references in relation to the construction of wavelets on general bounded domains for numerical DEs are \cite{ctu99, ds99, ks06}. In addition to the two important research directions described above, a substantial amount of noteworthy work has gone towards developing adaptive wavelet algorithms. Such a scheme is useful in handling problems with singularities and nonlinear problems. There is a subtle distinction between adaptive finite element and adaptive wavelet methods. The first one directly relies on iterative mesh refinements driven by a posteriori local errors; meanwhile the latter involves iterative trial space refinements by proper inclusions of extra basis functions \cite{d01}. Two of the seminal papers in this field are \cite{cdd00, cdd02}. Further developments of adaptive wavelet algorithms, to name a few, are reported in \cite{h99} for hyperbolic problems, \cite{dss09} for high dimensional elliptic problems, and \cite{ss09} for parabolic evolution problems. The use of wavelets in efficiently dealing with integral operators can be traced back to the influential papers \cite{abcr93, bcr91} and more recent results can be found in \cite{cf20}. Multilevel preconditioning techniques for linear systems stemming from Galerkin methods were studied in \cite{dk92}.

Towards the end of this paper, we shall consider the Helmholtz equation, which is a widely used model in wave propagation. This equation appears in seismic imaging (e.g., full waveform inversion), acoustics, and electromagnetism. One predominant challenge in solving the Helmholtz equation is caused by the infamous pollution effect; i.e., small enough $hk$ (here $h$ and $k$ refer to the mesh size and wave number respectively) provides no guarantee that the relative error has a bound independent of $k$ \cite{dms19}. To make matters worse, the linear system associated with the discretization process is ill-conditioned. There is an abundance of numerical schemes in the literature that deal with the Helmholtz equation with high wave numbers. For example, see \cite{gz19, hmp16, ms14} and references therein. To this day, developing a numerical scheme for the Helmholtz equation continues to be an active research area. In this paper, we aim to examine a wavelet-based method's potential in solving the one-dimensional Helmholtz equation with wave numbers of magnitude $\bo(10^4)$ or larger.

The paper is organized as follows. In Section~2, we shall review the main result of \cite{hm19}. In Section~3, we shall revisit the main result of \cite{hmw19} and provide a new example of a biorthogonal wavelet formed on the unit interval with homogeneous Dirichlet boundary conditions and maximum vanishing moments. In Section~4, we shall use this example to solve the Helmholtz equation with high wave numbers and consider the biharmonic equation. Finally, we outline some future work in Section~5.

\section{Derivative-orthogonal Riesz wavelets in $H^{m}(\R)$}
In this section, we review the main result presented in \cite{hm19}. For complete details and proofs, we refer interested readers to the foregoing paper.

To facilitate our discussion, we introduce some definitions. The Fourier transform of $f\in \Lp{1}$ is defined as $\wh{f}(\xi):=\int_{\R} f(x) e^{-ix\xi} dx$ for all $\xi \in \R$. We then extend the definition to tempered distributions by means of duality. Let $\tau \in \R$. The Sobolev space $H^{\tau}(\R)$ contains all tempered distributions on $\R$ such that
\[
\| f \|^2_{H^{\tau}(\R)} := \frac{1}{2 \pi} \int_{\R} |\wh{f}(\xi)|^{2}(1+|\xi|^{2})^{\tau} d\xi < \infty.
\]
Let $(l_{0}(\Z))^{r \times s}$ be the space of all finitely supported matrix-valued filters $a: \Z \rightarrow \C^{r \times s}$. For $a \in (l_{0}(\Z))^{r \times s}$, we define $\wh{a}(\xi):=\sum_{k \in \Z} a(k) e^{-ik\xi}$, which is an $r\times s$ matrix of $2\pi$-periodic trigonometric polynomials. A wavelet affine system and a Riesz wavelet in the Sobolev space $H^{\tau}(\R)$ are defined the same way as in Section~1; we only need to replace $m \in \N \cup \{0\}$ with $\tau \in \R$. A compactly supported Riesz wavelet $\{\phi;\psi\}$ is typically derived from a compactly supported refinable vector function via the following refinability structure:
\[
\wh{\phi}(2\xi) = \wh{a}(\xi) \wh{\phi}(\xi) \quad \text{and} \quad \wh{\psi}(2\xi) = \wh{b}(\xi) \wh{\phi}(\xi), \quad \text{a.e. } \xi \in \R
\]
for some filters $a \in (l_{0}(\Z))^{r\times r}$ and $b \in (l_{0}(\Z))^{s\times r}$.
For $J\in \Z$, we say that $(\AS_J^{-\tau}(\tilde{\phi};\tilde{\psi}), \AS_J^\tau(\phi;\psi))$ is a biorthogonal wavelet in $(H^{-\tau}(\R),H^{\tau}(\R))$ if (1) $\AS^{\tau}_{J}(\phi;\psi)$ is a Riesz basis for $H^{\tau}(\R)$ and $\AS^{-\tau}_{J}(\tilde{\phi};\tilde{\psi})$ is a Riesz basis for $H^{-\tau}(\R)$, as well as (2) $\AS^{\tau}_{J}(\phi;\psi)$ and $\AS^{-\tau}_{J}(\tilde{\phi};\tilde{\psi})$ are biorthogonal to each other.
It is known in \cite{han12,hanbook} that
$(\AS_J^{-\tau}(\tilde{\phi};\tilde{\psi}), \AS_J^\tau(\phi;\psi))$ is a biorthogonal wavelet in $(H^{-\tau}(\R),H^{\tau}(\R))$ for some $J\in \Z$ if and only if $(\AS_J^{-\tau}(\tilde{\phi};\tilde{\psi}), \AS_J^\tau(\phi;\psi))$ is a biorthogonal wavelet in $(H^{-\tau}(\R),H^{\tau}(\R))$ for all $J\in \Z$. Consequently,
we say that $(\{\tilde{\phi};\tilde{\psi}\},\{\phi;\psi\})$ is \emph{a biorthogonal wavelet in} $(H^{-\tau}(\R),H^{\tau}(\R))$ if $(\AS_0^{-\tau}(\tilde{\phi};\tilde{\psi}), \AS_0^\tau(\phi;\psi))$ is a biorthogonal wavelet in $(H^{-\tau}(\R),H^{\tau}(\R))$.

The bracket product is a key tool in the analysis of shift-invariant spaces.
Let $f \in (H^{\tau}(\R))^{r \times t}$ and $g \in (H^{-\tau}(\R))^{s \times t}$ be vector functions. For $\tau \in \R$, we define \emph{the bracket product} by
\[
[\wh{f},\wh{g}]_{\tau}(\xi):= \sum_{k \in \Z} \wh{f}(\xi + 2\pi k) \overline{\wh{g}(\xi + 2\pi k)}^{\tp}(1+|\xi|^2)^{\tau}, \quad \xi \in \R.
\]
When $\tau =0$, $[\wh{f},\wh{g}](\xi):=[\wh{f},\wh{g}]_{0}(\xi)$. We say that the integer shifts of a vector function $\phi\in (H^\tau(\R))^{r\times 1}$ are \emph{stable} in $H^\tau(\R)$ if there exists a positive constant $C$ such that
\[
C^{-1} I_r \le [\wh{\phi},\wh{\phi}]_\tau(\xi)\le C I_r, \qquad a.e.\, \xi\in \R.
\]
If $\phi\in H^\tau(\R)$ has compact support, then the integer shifts of $\phi$ are stable in $H^\tau(\R)$ if and only if
$\text{span}\{\wh{\phi}(\xi + 2\pi k) \setsp k\in \Z\}=\C^{r}$ for all $\xi \in \R$ (see \cite[Theorem~5.3.6]{hanbook}).
Next, we introduce the definitions of sum rules and vanishing moments. The former is closely related to polynomial reproduction orders, while the latter is linked to polynomial annihilation. For a given filter $a \in (l_{0}(\Z))^{r \times r}$, we say that $a$ satisfies \emph{order $m$ sum rules} if there is a matching filter $v \in  (l_{0}(\Z))^{1 \times r}$ such that $\wh{v}(0) \neq 0$,
\begin{equation} \label{sr}
\wh{v}(2\xi) \wh{a}(\xi) = \wh{v}(\xi) + \bo(|\xi|^{m}), \quad \wh{v}(2\xi)\wh{a}(\xi + \pi)=\bo(|\xi|^{m}), \quad \xi \rightarrow 0.
\end{equation}
We also use the notation $\sr(a)$ to denote the highest order of sum rules satisfied by a filter $a$. Here, $\wh{f}(\xi) = \wh{g}(\xi) + \bo(|\xi|^{m})$ as $\xi \rightarrow 0$ is equivalent to saying that $f^{(j)}(0) = g^{(j)}(0)$ for all $j=0,\dots,m-1$. On the other hand, given a compactly supported function $\psi \in (L_{2}(\R))^{s \times 1}$, we say that $\psi$ has \emph{order $n$ vanishing moments} if $\langle \psi,x^j \rangle = \int_{\R} \psi(x) x^j dx =0$ for all $j=0,\dots,n-1$; or in other words, $\wh{\psi}(0)=\dots=\wh{\psi}^{(n-1)}(0)=0$. Notation-wise, we write $\vmo(\psi)$ to indicate the highest vanishing moments attained by $\psi$.

Let $\phi$ be a tempered distribution on $\R$. Recall that the smoothness/regularity of $\phi$ is measured by its \emph{smoothness exponent}:
\[
\sm(\phi):=\sup\{\tau \in \R \setsp \phi \in H^{\tau}(\R)\}.
\]
By convention, we set $\sm(\phi):=-\infty$ if $\{\tau \in \R\setsp \phi \in H^{\tau}(\R)\} = \emptyset$.

The theorem immediately below provides a full characterization of Riesz wavelets in the Sobolev space $H^{\tau}(\R)$, where $\tau \in \R$. Riesz wavelets and dual wavelet frames in the Sobolev space $H^{\tau}(\R)$, where $\tau \in \R$, were initially studied in \cite{hs09}.

\begin{theorem} \label{thm:rieszsobolev}
	(a special case of \cite[Theorem 6.4.6]{hanbook})
	Let $a, b, \tilde{a}, \tilde{b} \in (l_{0}(\Z))^{r \times r}$.
Assume that
\begin{enumerate}
\item[(i)]
$1$ is a simple eigenvalue of $\wh{a}(0)$ and $\det(2^{j}I_{r}-\wh{a}(0))\neq 0$ for all $j\in \N$;
\item[(ii)] $1$ is a simple eigenvalue of $\wh{\tilde{a}}(0)$ and $\det(2^{j}I_{r}-\wh{\tilde{a}}(0))\neq 0$ for all $j\in \N$.
\end{enumerate}
Let $\phi, \tilde{\phi}$ be $r \times 1$ vectors of compactly supported distributions satisfying
\begin{equation} \label{I:phidphi}
\wh{\phi}(2\xi) = \wh{a}(\xi)\wh{\phi}(\xi), \quad \wh{\tilde{\phi}}(2\xi) = \wh{\tilde{a}}(\xi)\wh{\tilde{\phi}}(\xi).
\end{equation}
Define $\psi$ and $\tilde{\psi}$ by
	\begin{equation} \label{I:psidpsi}
	\wh{\psi}(\xi) = \wh{b}(\xi/2)\wh{\psi}(\xi/2), \quad \wh{\tilde{\psi}}(\xi) = \wh{\tilde{b}}(\xi/2)\wh{\tilde{\psi}}(\xi/2).
	\end{equation}
	For $\tau \in \R$, the pair $(\{\tilde{\phi};\tilde{\psi}\},\{\phi;\psi\})$ is a biorthogonal wavelet in $(H^{-\tau}(\R),H^{\tau}(\R))$ if and only if
\begin{itemize}
\item[(1)] $(\{\tilde{a};\tilde{b}\},\{a;b\})$ is a biorthogonal wavelet filter bank; that is,
	\[
	\left[ \begin{matrix} \wh{\tilde{a}}(\xi) &\wh{\tilde{a}}(\xi+\pi)\\
		\wh{\tilde{b}}(\xi) &\wh{\tilde{b}}(\xi+\pi)
	\end{matrix}\right]
	\left[ \begin{matrix}
		\ol{\wh{a}(\xi)}^\tp &\ol{\wh{b}(\xi)}^\tp\\
		\ol{\wh{a}(\xi+\pi)}^\tp &\ol{\wh{b}(\xi+\pi)}^\tp\end{matrix}\right]=
	I_{2r},
	\]
\item[(2)] $\ol{\wh{\phi}(0)}^\tp \wh{\tilde{\phi}}(0)=1$, $\phi \in (H^{\tau}(\R))^{r\times 1}$ and $\tilde{\phi} \in (H^{-\tau}(\R))^{r\times 1}$.
\item[(3)] The integer shifts of $\phi$ and $\tilde{\phi}$ are biorthogonal to each other:
\be \label{phi:bio}
\la \tilde{\phi}, \phi(\cdot-k)\ra:=
\int_{\R} \tilde{\phi}(x) \ol{\phi(x-k)}^\tp dx=\td(k) I_r,\qquad \forall\; k\in \Z,
\ee
where $\td(0)=1$ and $\td(k)=0$ for all $k\ne 0$.
\item[(4)] $\wh{\psi}(\xi) = o(|\xi|^{-\tau})$ as $\xi \rightarrow 0$ (i.e., $\vmo(\psi) > -\tau$) if $\tau \le 0$, and $\wh{\tilde{\psi}}(\xi) = o(|\xi|^{\tau})$ as $\xi \rightarrow 0$ (i.e., $\vmo(\tilde{\psi}) > \tau$) if $\tau \ge 0$.
\end{itemize}
\end{theorem}

The conditions in items (i) and (ii) of Theorem~\ref{thm:rieszsobolev} are not assumed in \cite[Theorem 6.4.6]{hanbook}. It is well known (e.g., see \cite[Theorem~5.1.3]{hanbook})
that the condition in item (i) guarantees the existence and uniqueness (up to a multiplicative constant) of a compactly supported refinable vector function $\phi$ satisfying $\wh{\phi}(2\xi)=\wh{a}(\xi)\wh{\phi}(\xi)$.
Item (4) puts conditions on vanishing moments. For $\tau>0$, no vanishing moments are required for $\psi$ at all, while for $\tau<0$,
no vanishing moments are required for
$\tilde{\psi}$.
If $\phi\in (H^\tau(\R))^{r\times 1}$ satisfies $\wh{\phi}(2\xi)=\wh{a}(\xi)\wh{\phi}(\xi)$ with $a\in (l_{0}(\Z))^{r \times r}$ and $\tau\ge 0$, then item (i) must hold if the integer shifts of $\phi$ are stable in $H^\tau(\R)$.
Note that the biorthogonality condition in \eqref{phi:bio} is equivalent to saying that $[\wh{\tilde{\phi}},\wh{\phi}](\xi)=I_r$ for almost every $\xi\in \R$.
Since $\phi\in (H^\tau(\R))^{r\times 1}$ and $\tilde{\phi}\in (H^{-\tau}(\R))^{r\times 1}$ have compact support, item (3) implies that the integer shifts of $\phi$ are stable in $H^\tau(\R)$, while the integer shifts of $\tilde{\phi}$ are stable in $H^{-\tau}(\R)$.
Moreover, the biorthogonality condition \eqref{phi:bio} in item (3) can be fully characterized by $\sm(a)>\tau$ and $\sm(\tilde{a})>-\tau$, where the smoothness exponent $\sm(a)$ is defined in \cite[(5.6.44)]{hanbook}.
For more details, see \cite[Theorem~6.4.5]{hanbook} and \cite{han01,han03jat}.

Built on Theorem~\ref{thm:rieszsobolev},
the following result
characterizes derivative-orthogonal Riesz wavelets.

\begin{theorem}
(\cite[Theorems~2 and~4]{hm19})
Let $\phi=(\phi_1,\ldots,\phi_r)^\tp$ be a compactly supported refinable vector function in $\HH{m}$ with $m\in \NN$ such that $\wh{\phi}(2\xi)=\wh{a}(\xi)\wh{\phi}(\xi)$ for some $a\in \lrs{0}{r}{r}$. Then
\begin{enumerate}
\item[(i)] there exists a finitely supported high-pass filter $b\in \lrs{0}{r}{r}$ such that $\{\phi;\psi\}$ with $\wh{\psi}(\xi):=\wh{b}(\xi/2)\wh{\phi}(\xi/2)$ is an $m$th-order derivative-orthogonal Riesz wavelet in the Sobolev space $\HH{m}$ satisfying \eqref{modo:cond:phi} and \eqref{modo:cond:psi} if and only if the integer shifts of $\phi$ are stable and the filter $a$ has at least order $2m$ sum rules (i.e., $\sr(a)\ge 2m$).
\item[(ii)] Under the condition that the integer shifts of $\phi$ are stable, for any $b\in \lrs{0}{r}{r}$, $\{\phi;\psi\}$ with $\wh{\psi}(\xi):=\wh{b}(\xi/2)\wh{\phi}(\xi/2)$ is an $m$th-order derivative-orthogonal Riesz wavelet in the Sobolev space $\HH{m}$  satisfying \eqref{modo:cond:phi} and \eqref{modo:cond:psi} if and only if
\[
\wh{b}(\xi) [\wh{\phi^{(m)}},\wh{\phi^{(m)}}](\xi)
\ol{\wh{a}(\xi)}^\tp+\wh{b}(\xi+\pi) [\wh{\phi^{(m)}},\wh{\phi^{(m)}}](\xi+\pi)
\ol{\wh{a}(\xi+\pi)}^\tp=0
\]
and
\[
\det(\{\wh{a};\wh{b}\})(\xi):=\det\left(\left[ \begin{matrix} \wh{a}(\xi) &\wh{a}(\xi+\pi)\\
\wh{b}(\xi) &\wh{b}(\xi+\pi)
\end{matrix}\right]\right)\ne 0, \qquad \forall\; \xi\in \R.
\]
are satisfied. Moreover, for every $J\in \Z$, $\AS^\tau_J(\phi;\psi)$ is a Riesz basis in the Sobolev space $\HH{\tau}$ for all $\tau$ in the nonempty open interval $(2m-\sm(\phi), \sm(\phi))$ with $\vmo(\psi)=\sr(a)-2m$.
\end{enumerate}
\end{theorem}

For scalar filters, there is in fact an explicit formula to find their corresponding $m$-th order derivative orthogonal Riesz wavelets. See \cite[Theorem 5]{hm19} for details. To make the presentation of this paper self-contained, we reproduce the second-order derivative-orthogonal Riesz wavelet generated from Hermite cubic splines, since it will be used in Section~4.2.

\begin{example} \label{ex:hmtcubdow}
	\normalfont
	Let $\phi:=(\phi_{1},\phi_{2})$ be the well-known Hermite cubic splines with $\phi_{1}$ and $\phi_{2}$ given below
	\[
	\phi_1(x)=\begin{cases}
	(1-x)^2(1+2x), &x\in [0,1],\\
	(1+x)^2(1-2x), &x\in [-1,0),\\
	0, &\text{otherwise},
	\end{cases}
	\qquad
	\phi_2(x)=\begin{cases}
	(1-x)^2x, &x\in [0,1],\\
	(1+x)^2x, &x\in [-1,0),\\
	0, &\text{otherwise}.
	\end{cases}
	\]
	The integer shifts of $\phi$ are indeed stable and $\sm(\phi)=2.5$. Moreover, $\sr(a)=4$, where the filter $a$ is defined below. Note that $\phi, \psi$ satisfy equations $\wh{\phi}(2\xi) = \wh{a}(\xi)\wh{\phi}(\xi)$ and $\wh{\psi}(\xi) = \wh{b}(\xi/2)\wh{\phi}(\xi/2)$ respectively with the filters $a,b \in (l_{0}(\Z))^{2 \times 2}$ defined as follows:
	\[
	a=\left\{ \begin{bmatrix} \tfrac{1}{4} &\tfrac{3}{8}\\[0.3em]
	-\tfrac{1}{16} &-\tfrac{1}{16}\end{bmatrix},\quad
	\begin{bmatrix} \tfrac{1}{2} &0 \\[0.3em]
	0 &\tfrac{1}{4}\end{bmatrix}, \quad
	\begin{bmatrix} \tfrac{1}{4} &-\tfrac{3}{8}\\[0.3em]
	\tfrac{1}{16} &-\tfrac{1}{16}\end{bmatrix}\right\}_{[-1,1]}, \quad
	b=\left\{\begin{bmatrix} \tfrac{1}{2} & 0 \\[0.3em]
	0 & \tfrac{1}{2}
	\end{bmatrix}\right\}_{[1,1]}.
	\]
	Then, $\{\phi;\psi\}$ is a second-order derivative-orthogonal Riesz wavelet in $H^{2}(\R)$. Also, the wavelet affine system $\AS_{0}^{\tau}(\phi;\psi)$ is a Riesz basis in $H^{\tau}(\R)$ for all $\tau \in (3/2,5/2)$. See Figure~\ref{fig:hmtcubdow} for plots of $\phi$ and $\psi$.
	
	\begin{figure}[htbp]
		\centering
		\begin{subfigure}[b]{0.28\textwidth}
			 \includegraphics[width=\textwidth]{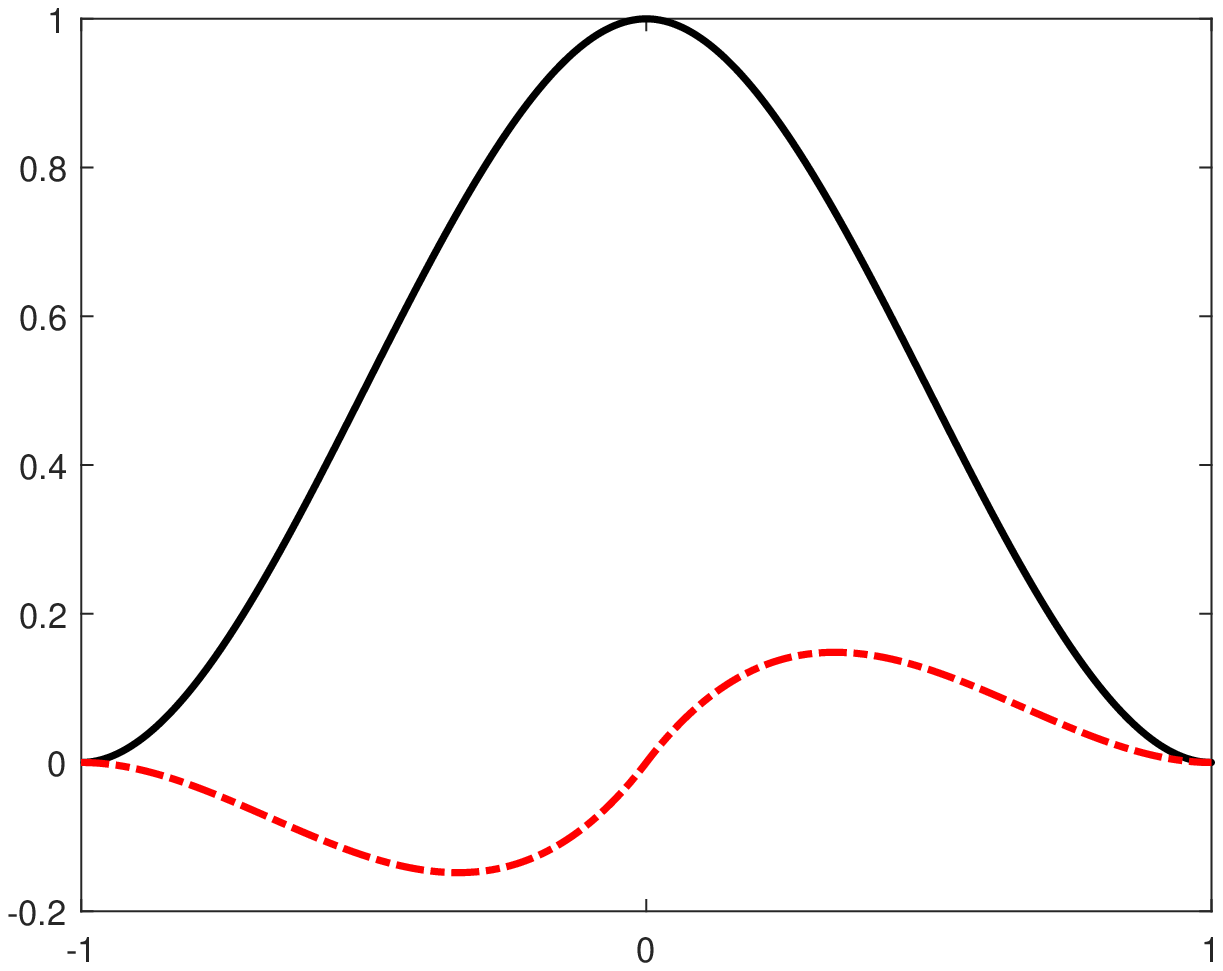}
		\end{subfigure}
		\begin{subfigure}[b]{0.28\textwidth}
			 \includegraphics[width=\textwidth]{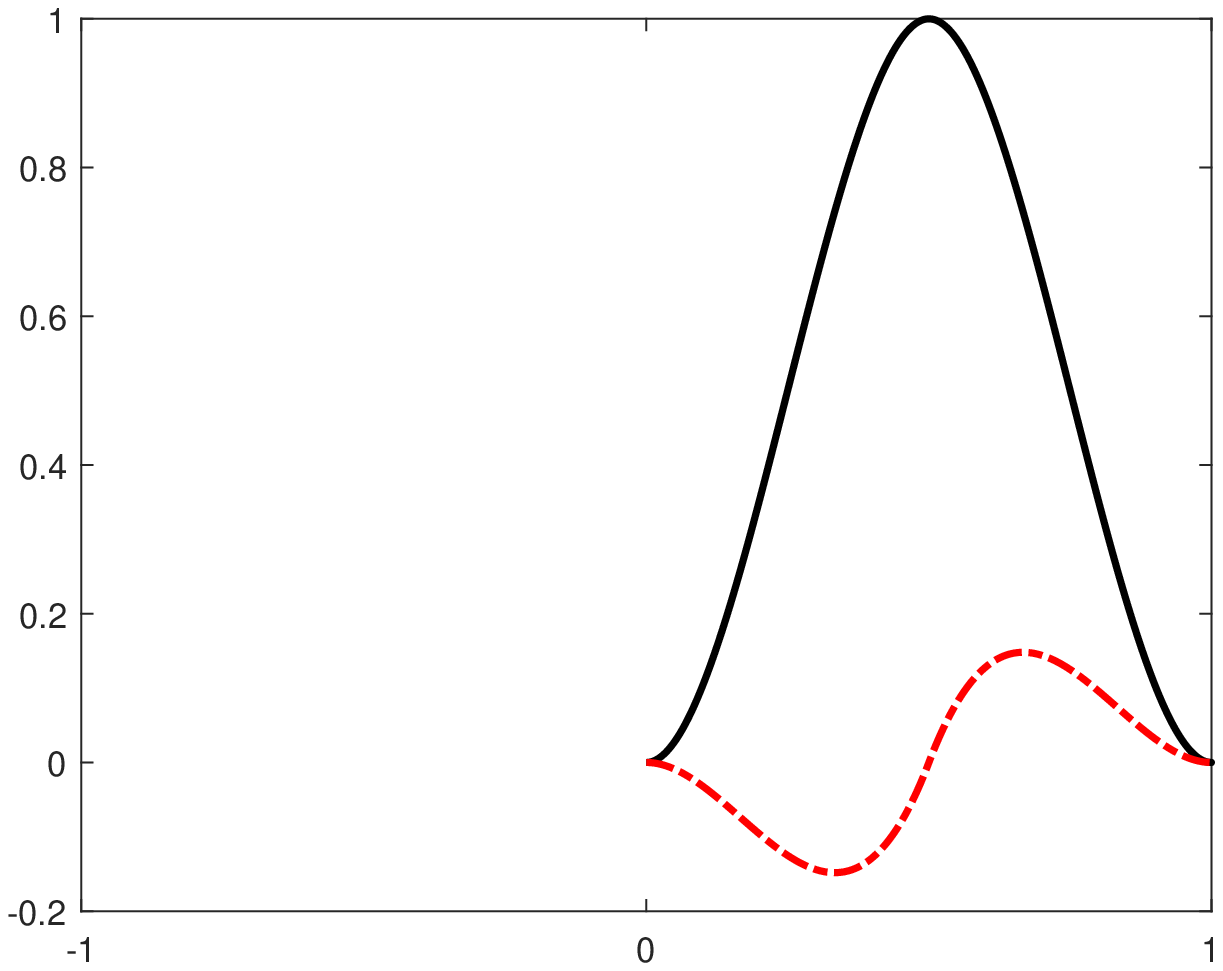}
		\end{subfigure}
		\caption{Black solid and red dashed lines correspond to the first and second components of a vector function respectively. Left: the refinable vector function $\phi=(\phi_{1},\phi_{2})$. Right: the wavelet vector function $\psi=(\psi_1, \psi_2)$ such that $\{\phi;\psi\}$ is a second-order derivative-orthogonal Riesz wavelet in $H^{2}(\R)$.}
		\label{fig:hmtcubdow}
	\end{figure}
\end{example}

Other examples of zeroth, first, and second-order derivative-orthogonal Riesz wavelets generated from B-splines and Hermite splines are presented in \cite[Section 3]{hm19}.

\section{Construction of biorthogonal wavelets in $L_{2}([0,N])$ with maximum vanishing moments and boundary conditions}

In this section, we shift our attention to a general construction of compactly supported biorthogonal wavelets in $L_{2}([0,N])$ with $N\in \N$. For complete details and proofs, we refer readers to \cite{hmw19}.

First and foremost, we review relevant studies that precede ours. A pioneering work in the construction of Daubechies orthogonal wavelets in $L_{2}([0,1])$ can be attributed to Meyer \cite{m91}, which was followed up by \cite{ahjp94, a93, cdv93, cq92, pst95} and other references. A more general construction of orthogonal wavelets in $L_{2}([0,N])$ is available in \cite{ak12}. The use of orthogonal wavelets in $L_2([0,N])$ within the framework of numerical DEs is not preferable for two reasons. Firstly, they often do not have any analytic expressions with a few exceptions like the Haar wavelet. Even then, this particular example has low smoothness. Secondly, their orthogonality property often conflicts with the boundary conditions. I.e., orthogonal wavelets often fail to satisfy the boundary conditions. If orthogonal wavelets in $L_2([0,N])$ would satisfy a given boundary condition, then they cannot have high orders of vanishing moments, which are the key property for sparse representations.
Consequently, implementing orthogonal wavelets in a numerical scheme becomes a cumbersome task. In contrast, biorthogonal wavelets are preferable mainly because they do not suffer from the previous disadvantages. Within the scope of biorthogonal wavelets, some existing constructions can solely be applied to a particular refinable function $\phi$; e.g., Hermite cubic splines in \cite{ahl17,cf16, cf17, dhjk00, ds10}. Scalar biorthogonal wavelets generated from B-splines can be formed on the unit interval by utilizing the methods in \cite{dku99, j09}. Additionally, the method explained in \cite{j09} yields spline wavelets on the unit interval with homogeneous boundary conditions of arbitrary order. On a similar note, \cite{ds98} provides a means to construct biorthogonal wavelets satisfying complementary boundary conditions. The constructions proposed in \cite{hanbook, hm18} rely on an operator introduced in \cite{cdv93}. They can be applied to a relatively large class of wavelets and framelets (i.e., a generalization of wavelets with redundancy) satisfying some symmetry property. However, these constructions sacrifice the vanishing moments near the boundary in order to fulfill the boundary conditions. Limited applicability as well as inability to concurrently retain maximum vanishing moments and meet the boundary conditions are two common deficiencies in the above constructions. This is where our contribution comes in.

To aid the discussion of our main result, we first introduce a notation and an important result that describes the structure of compactly supported Riesz wavelets on $[0,\infty)$. Suppose $f$ is a compactly supported vector function. Then, $\fs(f)$ is defined to be the shortest interval with integer endpoints such that $f$ vanishes outside $\fs(f)$. Also, denote $f_{j;k}:=2^{j/2}f(2^j \cdot -k)$.

The theorem below essentially guarantees the existence of compactly supported vector functions $\tilde{\phi}^{L}$, $\tilde{\psi}^{L}$ given compactly supported vector functions $\phi^{L}$, $\psi^{L}$: a critical fact for the construction procedure under consideration.

\begin{theorem}
(\cite[Theorem 2.2]{hmw19})
	Let $(\{\tilde{\phi};\tilde{\psi}\},\{\phi;\psi\})$ be a compactly supported biorthogonal wavelet in $\Lp{2}$ with $\phi,\psi,\tilde{\phi},\tilde{\psi}\in (\Lp{2})^r$. Then Theorem~\ref{thm:rieszsobolev} holds with $\tau =0$.
	Define
	\be \label{fsupp:phi}
	[l_\phi, h_\phi]:=\fs(\phi),\quad
	[l_\psi, h_\psi]:=\fs(\psi),\quad
	[l_{\tilde{\phi}},
	 h_{\tilde{\phi}}]:=\fs(\tilde{\phi}),\quad [l_{\tilde{\psi}}, h_{\tilde{\psi}}]:=\fs(\tilde{\psi})
	\ee
	and
	\be \label{fsupp:a}
	[l_a, h_a]:=\fs(a),\quad
	[l_b, h_b]:=\fs(b),\quad
	[l_{\tilde{a}}, h_{\tilde{a}}]:=\fs(\tilde{a}),\quad [l_{\tilde{b}}, h_{\tilde{b}}]:=\fs(\tilde{b}).
	\ee
	%
	Let $\phi^L$ and $\psi^L$ be vectors of compactly supported functions in $L_2([0,\infty))$. Define
	\be \label{PhiPsiI}
	\Phi:= \{\phi^{L}\} \cup \{\phi(\cdot -k): k \ge n_{\phi}\}, \quad
	\Psi:= \{\psi^{L}\} \cup \{\psi(\cdot -k): k \ge n_{\psi}\}
	\ee
	with $n_\phi\ge \max(-l_\phi,-l_a)$ and $n_\psi\ge \max(-l_\psi, \frac{n_\phi-l_b}{2})$.
Define
\[ \AS_J(\Phi;\Psi)_{[0,\infty)}:=\{2^{J/2}\varphi(2^J\cdot) \setsp
	\varphi \in \Phi\} \cup\{ 2^{j/2}\eta(2^j\cdot) \setsp j\ge J, \eta \in \Psi\},\qquad J\in \Z.
	\]
	Suppose that $\AS_0(\Phi;\Psi)_{[0,\infty)}$ is a Riesz basis of $L_2([0,\infty))$
	and satisfies
	\begin{align}
	 &\phi^L=2A_L\phi^L(2\cdot)+2\sum_{k=n_\phi}^{m_\phi} A(k) \phi(2\cdot-k),
	\nonumber\\
	 &\psi^L=2B_L\phi^L(2\cdot)+2\sum_{k=n_\phi}^{m_\psi} B(k) \phi(2\cdot-k),
	\label{I:psi}
	\end{align}
	for some matrices $A_L,B_L$ and finitely supported sequences $A,B$ of matrices.
	Then
	\begin{enumerate}
		\item[(1)] there must exist compactly supported vector functions $\tilde{\phi}^L, \tilde{\psi}^L$ in $L_2([0,\infty))$ and integers $n_{\tilde{\phi}}\ge \max(-l_{\tilde{\phi}},-l_{\tilde{a}},n_\phi)$ and $n_{\tilde{\psi}}\ge \max(-l_{\tilde{\psi}},\frac{n_{\tilde{\phi}}-l_{\tilde{b}}}{2},n_\psi)$ such that
		 $\AS_0(\tilde{\Phi};\tilde{\Psi})_{[0,\infty)}$ is the dual Riesz basis of $\AS_0(\Phi;\Psi)_{[0,\infty)}$ in $L_2([0,\infty))$,
		where $\AS_0(\tilde{\Phi};\tilde{\Psi})_{[0,\infty)}:=\tilde{\Phi} \cup \{ 2^{j/2} \tilde{\eta}(2^j\cdot) \setsp j\in \NN, \tilde{\eta}\in \tilde{\Psi}\}$ with $\NN:=\N\cup\{0\}$ and
		\be \label{tPhiPsiI}
		\tilde{\Phi}:=\{ \tilde{\phi}^L\}\cup\{ \tilde{\phi}(\cdot-k) \setsp k\ge n_{\tilde{\phi}}\},\quad
		\tilde{\Psi}:=\{ \tilde{\psi}^L\}\cup\{ \tilde{\psi}(\cdot-k) \setsp k\ge n_{\tilde{\psi}}\};
		\ee

		\item[(2)] there exist matrices $\tilde{A}_L,\tilde{B}_L$ and finitely supported sequences $\tilde{A},\tilde{B}$ of matrices such that
		\begin{align}
		 &\tilde{\phi}^L=2\tilde{A}_L\tilde{\phi}^L(2\cdot)+2\sum_{k=n_{\tilde{\phi}}}^{m_{\tilde{\phi}}} \tilde{A}(k) \tilde{\phi}(2\cdot-k), \nonumber\\
		 &\tilde{\psi}^L=2\tilde{B}_L\tilde{\phi}^L(2\cdot)+2\sum_{k=n_{\tilde{\phi}}}^{m_{\tilde{\psi}}} \tilde{B}(k) \tilde{\phi}(2\cdot-k),
		\label{I:psi:dual}
		\end{align}
		and
		\begin{align*}
		 &\tilde{\phi}(\cdot-k_0)=2\sum_{k=n_{\tilde{\phi}}}^\infty
		\tilde{a}(k-2k_0) \tilde{\phi}(2\cdot-k),\qquad \forall\; k_0\ge n_{\tilde{\phi}},\\
		 &\tilde{\psi}(\cdot-k_0)=2\sum_{k=n_{\tilde{\phi}}}^\infty
		\tilde{b}(k-2k_0) \tilde{\phi}(2\cdot-k),\qquad \forall\; k_0\ge n_{\tilde{\psi}}.
		\end{align*}
	\end{enumerate}
\end{theorem}

Before we present the complete algorithms, we shall outline several key steps to better illustrate how the construction operates. For the following, let $(\{\tilde{\phi};\tilde{\psi}\},\{\phi;\psi\})$ be a compactly supported biorthogonal wavelet in $L_{2}(\R)$.
\begin{itemize}
	\item[(1)] We commence by forming a compactly supported $L_{2}([0,\infty))$ biorthogonal wavelet.
	\begin{itemize}
		\item[(1.1)] \textbf{(Algorithm~\ref{alg:Phi})} Construct the space $\Phi$ as in \eqref{PhiPsiI} satisfying a prescribed polynomial reproduction property (cannot be higher than that of $\phi$) and the refinability condition.
		\item[(1.2)] \textbf{(Algorithm~\ref{alg:integral})} Compute the inner products between the shifted versions of $\phi$ and $\tilde{\phi}$ restricted to $[0,1]$. This will help us to recover any inner products between the shifts of $\phi$ and $\tilde{\phi}$.
		\item[(1.3)] \textbf{(Algorithm~\ref{alg:Phidual})} Construct the space $\tilde{\Phi}$ as in \eqref{tPhiPsiI} satisfying a prescribed polynomial reproduction (cannot be higher than that of $\tilde{\phi}$), the refinability condition, and the biorthogonality condition with respect to $\Phi$.
		\item[(1.4)] \textbf{(Algorithm~\ref{alg:Psi})} Construct the space $\Psi$ as in \eqref{PhiPsiI} such that all elements in $\Psi$ are perpendicular to all elements in $\tilde{\Phi}$. Additionally, we need to ensure that $\psi^{L}$ and all interior primal wavelets, $\{\psi(\cdot -k): k \ge n_{\psi}\}$, are linearly independent. Construct the space $\tilde{\Psi}$ as in \eqref{tPhiPsiI} such that the biorthogonality condition with respect to $\Psi$ is satisfied and all elements in $\tilde{\Psi}$ are perpendicular to all elements in $\Phi$. Additionally, we need to ensure that $\tilde{\psi}^{L}$ and all interior dual wavelets, $\{\tilde{\psi}(\cdot -k): k \ge n_{\tilde{\psi}}\}$, are linearly independent.
	\end{itemize}
	\item[(2)] Repeat Step (1) for the reflected compactly supported $L_{2}(\R)$ biorthogonal wavelet $(\{\tilde{\mathring{\phi}},\tilde{\mathring{\psi}}\},\{\mathring{\phi},\mathring{\psi}\})$, where $\mathring{\phi}:=\phi(-\cdot)$, $\mathring{\psi}:=\psi(-\cdot)$, $\tilde{\mathring{\phi}}:=\tilde{\phi}(-\cdot)$, and  $\tilde{\mathring{\psi}}:=\tilde{\psi}(-\cdot)$, in order to obtain compactly supported $\mathring{\Phi}$, $\mathring{\Psi}$, $\tilde{\mathring{\Phi}}$, and $\tilde{\mathring{\Psi}}$.
	\item[(3)] Set $J_{0}$ to be the smallest nonnegative integer such that for all $j \ge J_{0}$,
	\[
	\max(m_\phi+n_{\mathring{\phi}}, m_\psi+n_{\mathring{\phi}}, m_{\mathring{\phi}}+n_\phi, m_{\mathring{\psi}}+n_\phi)\le 2^{j+1}N
	\]
	and each element in $\Phi(2^j)$, $\Psi(2^j)$, $\tilde{\Phi}(2^j)$, $\tilde{\Psi}(2^j)$ does not essentially touch both endpoints $0$ and $N$ simultaneously. I.e., neither $0$ nor $N$ is an interior point of $\fs(h)$ for all $h \in \Phi(2^j \cdot), \Psi(2^j \cdot), \tilde{\Phi}(2^j \cdot), \tilde{\Psi}(2^j \cdot)$. Similarly, set $\tilde{J}_{0}$ to be the smallest nonnegative integer such that for all $j \ge \tilde{J}_{0}$
	\[
	 \max(m_{\tilde{\phi}}+n_{\tilde{\mathring{\phi}}},
	 m_{\tilde{\psi}}+n_{\tilde{\mathring{\phi}}},
	 m_{\tilde{\mathring{\phi}}}+n_{\tilde{\phi}},
	 m_{\tilde{\mathring{\psi}}}+n_{\tilde{\phi}})
	\le 2^{j+1}N,
	\]
	each element in $\mathring{\Phi}(2^j)$, $\mathring{\Psi}(2^j)$, $\tilde{\mathring{\Phi}}(2^j)$, $\tilde{\mathring{\Psi}}(2^j)$ does not essentially touch both endpoints $0$ and $N$ simultaneously, as well as
	\[
	\begin{split}
		\la \psi^L_{j;0}, \tilde{\phi}^R_{j;2^jN-N}\ra=0,\quad
		\la \psi^L_{j;0}, \tilde{\psi}^R_{j;2^jN-N}\ra=0,
		\quad
		\la \tilde{\psi}^L_{j;0}, \phi^R_{j;2^jN-N}\ra=0,\quad
		\la \tilde{\psi}^L_{j;0}, \psi^R_{j;2^jN-N}\ra=0,\\
		\la \tilde{\phi}^L_{j;0}, \psi^R_{j;2^jN-N}\ra=0,\quad
		\la \tilde{\psi}^L_{j;0}, \psi^R_{j;2^jN-N}\ra=0,
		\quad
		\la \phi^L_{j;0}, \tilde{\psi}^R_{j;2^jN-N}\ra=0,\quad
		\la \psi^L_{j;0}, \tilde{\psi}^R_{j;2^jN-N}\ra=0,
	\end{split}
	\]
	where
	\[
	 \phi^R:=\mathring{\phi}^L(N-\cdot),\quad
	 \psi^R:=\mathring{\psi}^L(N-\cdot),\quad
	 \tilde{\phi}^R:=\tilde{\mathring{\phi}}^L(N-\cdot),\quad
	 \tilde{\psi}^R:=\tilde{\mathring{\psi}}^L(N-\cdot).
	\]
\end{itemize}
By \cite[Theorem 2.5]{hmw19}, the above Steps (1) and (2) ensure that $\AS_{J}(\tilde{\Phi},\tilde{\Psi})|_{[0,\infty)}$ and $\AS_{J}(\Phi,\Psi)|_{[0,\infty)}$ form a pair of compactly supported Riesz bases in $L_{2}([0,\infty))$. Similarly, $\AS_{J}(\tilde{\mathring{\Phi}},\tilde{\mathring{\Psi}})|_{[0,\infty)}$ and $\AS_{J}(\mathring{\Phi},\mathring{\Psi})|_{[0,\infty)}$ form a pair of compactly supported Riesz bases in $L_{2}([0,\infty))$. Without loss of generality, suppose that $\tilde{J}_{0} \ge J_{0}$. Define
\begin{align}
&\Phi_j:=\{\phi^L_{j;0}\} \cup \{ \phi^R_{j;2^jN-N} \}\cup \{\phi_{j;k} \setsp n_\phi\le k\le 2^jN-n_{\mathring{\phi}}\}, \label{Phij}\\
&\Psi_j:=\{ \psi^L_{j;0} \} \cup \{ \psi^R_{j;2^jN-N}\}\cup \{\psi_{j;k} \setsp n_\psi\le k\le 2^jN-n_{\mathring{\psi}}\}, \label{Psij}\\
&\tilde{\Phi}_j:=\{\tilde{\phi}^L_{j;0}\} \cup \{ \tilde{\phi}^R_{j;2^jN-N}\}\cup \{\tilde{\phi}_{j;k} \setsp n_{\tilde{\phi}}\le k\le 2^jN-n_{\tilde{\mathring{\phi}}}\}, \label{tPhij}\\
&\tilde{\Psi}_j:= \{ \tilde{\psi}^L_{j;0}\} \cup \{ \tilde{\psi}^R_{j;2^jN-N} \}\cup \{ \tilde{\psi}_{j;k} \setsp n_{\tilde{\psi}}\le k\le 2^jN-n_{\tilde{\mathring{\psi}}}\} \label{tPsij}.
\end{align}
Furthermore, define $\cB_J:=\Phi_J \cup\{ \Psi_j \setsp\; j\ge J\}$, and $\tilde{\cB}_J:=\tilde{\Phi}_J \cup\{ \tilde{\Psi}_j \setsp\; j\ge J\}$. Following the above three key steps, we have by \cite[Theorem 4.1]{hmw19} that $(\tilde{\cB}_{J},\cB_{J})$ constitutes a pair of biorthogonal Riesz bases of $L_{2}([0,N])$ for all $J \ge \tilde{J}_{0}$. Moreover, there exist matrices $A_{j}$, $B_{j}$, $\tilde{A}_{j}$, $\tilde{B}_{j}$ such that the following refinable structures hold
\[
\Phi_{j} = A_{j} \Phi_{j+1}, \quad \Psi_{j} = B_{j} \Phi_{j+1}, \quad \tilde{\Phi}_{j} = \tilde{A}_{j} \tilde{\Phi}_{j+1}, \quad \tilde{\Psi}_{j} = \tilde{B}_{j} \tilde{\Phi}_{j+1},
\]
and $[\overline{A_{j}}^{\tp},\overline{B_{j}}^{\tp}]$ is an invertible square matrix with
\[
\left[ \begin{matrix} \tilde{A}_j\\ \tilde{B}_j\end{matrix}\right]=[\ol{A_j}^\tp, \ol{B_j}^\tp]^{-1}
\]
for $j \ge J$. Note that orthogonal wavelets can also be constructed by the foregoing procedure. We now present Algorithms~\ref{alg:Phi} to \ref{alg:Psi} in their full form.

\begin{algorithm}
	(\cite[Algorithm 1]{hmw19}) Let $\phi\in (\Lp{2})^r$ be a compactly supported refinable vector function such that $\phi=2\sum_{k\in \Z} a(k) \phi(2\cdot-k)$ for some finitely supported filter $a\in \lrs{0}{r}{r}$ and $a$ has $m$ sum rules in \eqref{sr}
	with respect to a moment matching filter $\vgu\in \lrs{0}{1}{r}$ satisfying $\wh{\vgu}(0)\wh{\phi}(0)=1$.
	Define $[l_\phi,h_\phi]:=\fs(\phi)$ and $[l_a,h_a]:=\fs(a)$.
	\begin{enumerate} \label{alg:Phi}
		\item[(S1)] Choose $n_\phi\ge \max(-l_\phi,-l_a)$. We often set $n_\phi:=\max(-l_\phi,-l_a)$.
		\item[(S2)] Define $\phi^c$ to be the column vector function consisting of $\phi(\cdot-k)\chi_{[0,\infty)}$ for $k$ decreasing from $n_\phi-1$ to $1-h_\phi$.
		Then
		\be \label{phic:refstr}
		\phi^c=2E_c \phi^c(2\cdot)+2\sum_{k=n_\phi}^\infty E(k) \phi(2\cdot-k),
		\ee
		where $E_c=(a(k-2n))_{n_\phi-1\ge n,k\ge 1-h_\phi}$ and $E(k):=(a(k-2n))_{n_\phi-1\ge n\ge 1-h_\phi}$ for $k \ge n_\phi$. If all the entries in $\phi^c$ are not linearly independent, then we delete as many entries as possible from $\phi^c$ so that all the deleted entries are linear combinations of entries kept. The relation \eqref{phic:refstr} still holds after appropriate modification.
		
		\item[(S3)]
		Let $\pp(x):=[x^{j_0},\ldots,x^{j_n}]^\tp$ with $j_0,\ldots,j_n\in \{0,\ldots,m-1\}$ (To preserve polynomial reproduction property, we often take $\pp(x)=[1,x,\ldots,x^{m-1}]^\tp$).
		Define a matrix $A_{\pp}$ via
		\[
		A_{\pp} \phi^c:=
		 \sum_{k=1-h_{\phi}}^{n_{\phi}-1}
		\sum_{j=0}^{m-1}
		\frac{(-i)^j}{j!}\pp^{(j)}(k) \wh{\vgu}^{(j)}(0)\phi(\cdot-k)\chi_{[0,\infty)}.
		\]
		%
		Perform row operations on $A_{\pp}$ to reduce it into row echelon form $A_r$.
		Define a (column) vector function
		$\phi^L:=A_c\phi^c$, where $A_c$ is an undetermined matrix in row echelon form with all leading coefficients being $1$ such that its first $m$ rows are given by $A_r$.
		
		\item[(S4)] Obtain a unique matrix $A_L^\tp$ through column operations by using the leading coefficient $1$ in $A_c^\tp$ to eliminate all other nonzero entries in $E_c^\tp A_c^\tp$ at the same row.
		Determine parameters in $A_c$ by solving
		$A_c E_c=A_L A_c$. In particular, if we take the particular choice
		\[
		\phi^L=A_{\pp}\phi^c \quad\mbox{with}\quad
		 \pp(x):=[x^{j_0},\ldots,x^{j_n}]^\tp\quad \mbox{and}\quad j_0,\ldots, j_n\in \{0,\ldots,m-1\},
		\]
		then $A_c E_c=A_L A_c$ automatically holds with $A_c:=A_{\pp}$ and $A_L:=\mbox{diag}(2^{-j_0},\ldots,2^{-j_n})$.
	\end{enumerate}
\end{algorithm}

\begin{algorithm} \label{alg:integral}
	(\cite[Theorem 3.2]{hmw19}) Let $\phi,\tilde{\phi}$ be two $r\times 1$ vectors of compactly supported functions in $\Lp{2}$ such that $\phi=2\sum_{k\in\Z} a(k) \phi(2\cdot-k)$ and $\tilde{\phi}=2\sum_{k\in\Z} \tilde{a}(k) \tilde{\phi}(2\cdot-k)$ for some finitely supported filters $a,\tilde{a}\in \lrs{0}{r}{r}$. Assume that $\wh{\phi}(0)\ne 0$ and $\wh{\tilde{\phi}}(0)\ne 0$.
	Define $[l_\phi,h_\phi]:=\fs(\phi)$ and $[l_{\tilde{\phi}},h_{\tilde{\phi}}]:=\fs(\tilde{\phi})$.
	\begin{enumerate}
		\item[(S1)] Define two vector functions by $\vec{\phi}:=[\phi(\cdot-1+h_\phi)\chi_{[0,1]},\ldots, \phi(\cdot+l_\phi)\chi_{[0,1]}]^\tp$ and $\vec{\tilde{\phi}}:=[\tilde{\phi}(\cdot-1+h_{\tilde{\phi}})
		\chi_{[0,1]},\ldots, \tilde{\phi}(\cdot+l_{\tilde{\phi}})\chi_{[0,1]}]^\tp$. Then
		\be \label{vecphi:refeq}
		\vec{\phi}=2A_0 \vec{\phi}(2\cdot)+2A_1 \vec{\phi}(2\cdot-1) \quad \mbox{and}\quad \vec{\tilde{\phi}}=2\tilde{A}_0 \vec{\tilde{\phi}}(2\cdot)+2\tilde{A}_1 \vec{\tilde{\phi}}(2\cdot-1)
		\ee
		with $A_\gamma:=(a(k+\gamma-2j))_{1-h_\phi\le j, k\le -l_\phi}$ and $\tilde{A}_\gamma:=(\tilde{a}(k+\gamma-2j))_{1-h_{\tilde\phi}\le j, k\le -l_{\tilde{\phi}}}$ for $\gamma=0,1$.
		\item[(S2)] If all the entries in $\vec{\phi}$ are not linearly independent on $[0,1]$, then we delete as many entries as possible from $\vec{\phi}$ so that all the deleted entries are linear combinations of entries kept. Do the same for $\vec{\tilde{\phi}}$. Then \eqref{vecphi:refeq} still holds with $A_0, A_1, \tilde{A}_0$ and $\tilde{A}_1$ being appropriately modified.
		
		\item[(S3)] Define $M:=\la \vec{\phi}, \vec{\tilde{\phi}} \ra:=\int_0^1 \vec{\phi}(x) \ol{\vec{\tilde{\phi}}(x)}^\tp dx$. Then the matrix $M$ is uniquely determined by the system of linear equations given by
		\[
		M=2A_0 M \ol{\tilde{A}_0}^\tp+2A_1 M \ol{\tilde{A}_1}^\tp
		\]
		under the normalization condition
		\[
		\vec{v} M \ol{\vec{\tilde{v}}}^\tp=1,
		\]
		where $\vec{v}$ is the unique row vector satisfying $\vec{v}(A_0+A_1)=\vec{v}$ and $\vec{v}\wh{\vec{\phi}}(0)=1$, while similarly $\vec{\tilde{v}}$ is the unique row vector satisfying $\vec{\tilde{v}}(\tilde{A}_0+\tilde{A}_1)=\vec{\tilde{v}}$ and $\vec{\tilde{v}}\wh{\vec{\tilde{\phi}}}(0)=1$.
	\end{enumerate}
\end{algorithm}

\begin{algorithm} \label{alg:Phidual}
	(\cite[Algorithm 3]{hmw19}) Let $(\{\tilde{\phi};\tilde{\psi}\},\{\phi;\psi\})$ be a compactly supported biorthogonal wavelet in $\Lp{2}$ associated with a finitely supported biorthogonal wavelet filter bank $(\{\tilde{a};\tilde{b}\},\{a;b\})$.
	Let $0\le m\le \sr(a)$ and $0\le \tilde{m}\le \sr(\tilde{a})$.
	Assume that
	$\Phi=\{\phi^L\}\cup\{\phi(\cdot-k) \setsp k\ge n_\phi\}$ is constructed by Algorithm~\ref{alg:Phi}.
	Define $[l_{\tilde{\phi}},h_{\tilde{\phi}}]:=\fs(\tilde{\phi})$ and $[l_{\tilde{a}},h_{\tilde{a}}]:=\fs(\tilde{a})$.
	\begin{enumerate}
		\item[(S1)] Choose $n_{\tilde{\phi}}\ge \max(-l_{\tilde{\phi}},-l_{\tilde{a}},n_\phi)$ such that $n_{\tilde{\phi}}$ is the smallest integer satisfying $\la \tilde{\phi}(\cdot-k), \phi^L\ra=0$ for all $k\ge n_{\tilde{\phi}}$.
		
		\item[(S2)] Define $\tilde{\phi}^c$ to be the vector function consisting of $\tilde{\phi}(\cdot-k)\chi_{[0,\infty)}$ for $k$ decreasing from $n_{\tilde{\phi}}-1$ to $1-h_{\tilde{\phi}}$. Then
		\be \label{tphic:refstr}
		\tilde{\phi}^c=2\tilde{E}_c \tilde{\phi}^c(2\cdot)+2\sum_{k=n_{\tilde{\phi}}}^\infty \tilde{E}(k) \tilde{\phi}(2\cdot-k),
		\ee
		where $\tilde{E}_c=(\tilde{a}(k-2n))_{n_{\tilde{\phi}}-1 \ge n,k\ge
			1-h_{\tilde{\phi}}}$ and $\tilde{E}(k):=(\tilde{a}(k-2n))_{n_{\tilde{\phi}}-1\ge n \ge 1-h_{\tilde{\phi}}}$ for $k \ge n_{\tilde{\phi}}$. If all the entries in $\tilde{\phi}^c$ are not linearly independent, then we delete as many entries as possible from $\tilde{\phi}^c$ so that all the deleted entries are linear combinations of entries kept. The relation \eqref{tphic:refstr} still holds after appropriate modification.
		
		\item[(S3)] Since $n_{\tilde{\phi}}\ge n_\phi$, we define a vector function $\mathring{\phi}^L$ by appending $\phi^L$ with $\phi(\cdot-k), n_\phi\le k< n_{\tilde{\phi}}$. Use Algorithm~\ref{alg:integral} to calculate $\la \tilde{\phi}^c, \mathring{\phi}^L\ra$.
		Define a vector function $\tilde{\phi}^L:=\tilde{A}_c \tilde{\phi}^c$ with $\#\tilde{\phi}^L=\#\mathring{\phi}^L$, where  the unknown $(\#\mathring{\phi}^L)\times (\#\tilde{\phi}^c)$ matrix $\tilde{A}_c$ is determined by solving the system of linear equations: $\tilde{A}_c \la \tilde{\phi}^c, \mathring{\phi}^L\ra=I_{\#\mathring{\phi}^L}$.
		
		\item[(S4)] Let $\pp(x):=[x^{\tilde{j}_0},\ldots, x^{\tilde{j}_{\tilde{n}}}]^\tp$ with $\tilde{j}_0,\ldots,\tilde{j}_{\tilde{n}}\in\{0,\ldots,\tilde{m}-1\}$ (To preserve polynomial reproduction property, we often take $\pp(x)=[1,x,\ldots,x^{\tilde{m}-1}]^\tp$).
		Find a matrix $\tilde{A}_p$ such that
		\[
		\tilde{A}_p \tilde{\phi}^c:=
		 \sum_{k=1-h_{\tilde{\phi}}}^{n_{\tilde{\phi}}-1} \sum_{j=0}^{\tilde{m}-1}
		\frac{(-1)^j}{j!}\pp^{(j)}(k) \wh{\tilde{\vgu}}^{(j)}(0)\tilde{\phi}(\cdot-k)\chi_{[0,\infty)},
		\]
		where $\tilde{\vgu}\in \lrs{0}{1}{r}$ is the moment matching filter for the sum rules of the filter $\tilde{a}$.
		Solve the linear equations
		$\tilde{A}_p=\la \pp, \mathring{\phi}^L\ra \tilde{A}_c$
		to further reduce the free parameters in $\tilde{A}_c$.
		
		\item[(S5)] Solve the equations $\tilde{A}_c\tilde{E}_c=\tilde{A}_c \tilde{E}_c \la \tilde{\phi}^c,\mathring{\phi}^L\ra \tilde{A}_c$ for the rest of free parameters in $\tilde{A}_c$.
	\end{enumerate}
\end{algorithm}

\begin{algorithm} \label{alg:Psi}
(\cite[Theorem 2.5]{hmw19})
Let $(\{\tilde{\phi};\tilde{\psi}\},\{\phi;\psi\})$ be a compactly supported biorthogonal wavelet in $\Lp{2}$ with $\phi,\psi,\tilde{\phi},\tilde{\psi}\in (\Lp{2})^r$. Define $l_\phi,l_\psi,l_{\tilde{\phi}}, l_{\tilde{\psi}}$ as in \eqref{fsupp:phi} and $l_a,l_b,l_{\tilde{a}}, l_{\tilde{b}}$ as in \eqref{fsupp:a}.
	\begin{enumerate}
		\item[(S1)]
		Define $n_\psi:=\max(-l_\psi,\lceil \frac{n_\phi-l_b}{2}\rceil)$ and $m_\phi:=k_\phi+\max(h_a-l_{\tilde{a}},0)$ with $k_\phi:=\max(2n_\phi+h_{\tilde{a}},2n_\psi+h_{\tilde{b}})-1$.
		Let $\eta$ be a column vector formed by listing all the entries in $\phi^L(2\cdot)$ and $\phi(2\cdot-k), k=n_\phi,\ldots,m_\phi$.
		Let $X$ be a matrix whose rows form a basis for the linear space of all row vectors $c$ satisfying
		\[
		\la c\eta, \tilde{\phi}^L\ra=0 \quad \mbox{and}\quad
		\la c\eta, \tilde{\phi}(\cdot-k)\ra=0,\qquad k=n_{\tilde{\phi}},\ldots,m_\phi+h_\phi-l_{\tilde{\phi}}-1.
		\]
		Let $Y$ be the matrix consisting of all row vectors $c$ such that
		$c \eta$ agrees with some entry of $\psi(\cdot-k), k=n_\psi,\ldots,\lceil \frac{m_\phi-l_b}{2}\rceil$. Write $Y=UX$ for some matrix $U$. Choose a matrix $V$ such that the square matrix $\left[ \begin{matrix} U\\ V\end{matrix}\right]$ is invertible.
		Define  $\Psi:=\{\psi^L\} \cup \{\psi(\cdot-k) \setsp k\ge n_{\psi}\}$, where $\psi^L:=VX \eta$ and can be rewritten in the form of \eqref{I:psi}.
		The matrix $V$ is often chosen so that $\psi^L$ has short support, satisfies some boundary conditions, or has a small condition number for
		the Riesz sequence $\{\psi^L\}\cup\{\psi(\cdot-k) \setsp k\ge n_\psi\}$.
		\item[(S2))] Define $n_{\tilde{\psi}}\ge \max(-l_{\tilde{\psi}},\lceil \frac{n_{\tilde{\phi}}
			-l_{\tilde{b}}}{2}\rceil, n_\psi)$ and $m_{\tilde{\phi}}:=\max(2n_{\tilde{\phi}}
		 +h_{a},2n_{\tilde{\psi}}+h_{b})+\max(h_{\tilde{a}}-l_a,0)-1$.
		Let $\tilde{\eta}$ be a column vector formed by listing all the entries in $\tilde{\phi}^L(2\cdot)$ and $\tilde{\phi}(2\cdot-k), k=n_{\tilde{\phi}},\ldots,m_{\tilde{\phi}}$.
		For each element in $h\in \{\psi^L\}\cup\{\psi(\cdot-k) \setsp k=n_\psi,\ldots,n_{\tilde{\psi}}-1\}$, there exists a unique element $\tilde{h}$ such that $\tilde{h}=\tilde{c} \tilde{\eta}$ with the coefficient row vector $\tilde{c}$ being uniquely determined by
		\[
		\la \tilde{c}\tilde{\eta}, h\ra=1\quad \mbox{and}\quad
		\la \tilde{c}\tilde{\eta}, g\ra=0, \qquad \forall\;
		g\in (\Phi \cup \Psi)\bs \{h\}.
		\]
		Define $\tilde{\Psi}:=\{\tilde{\psi}^L\}\cup \{ \tilde{\psi}(\cdot-k) \setsp k\ge n_{\tilde{\psi}}\}$, where $\tilde{\psi}^L$ is the vector of all above such $\tilde{h}$ and can be rewritten in the form of \eqref{I:psi:dual}.
	\end{enumerate}
\end{algorithm}

Various examples of compactly supported biorthogonal wavelets in $L_{2}([0,1])$ are available in \cite[Section 5]{hmw19}. We next turn to a compactly supported biorthogonal scalar wavelet on $L_{2}([0,1])$ satisfying homogeneous Dirichlet boundary conditions and maximum vanishing moments. Such an example will be used later in our numerical experiments for solving Helmholtz equations.

\begin{example} \label{ex:legall}
	\normalfont Consider the scalar biorthogonal wavelet $(\{\tilde{\phi};\tilde{\psi}\},\{\phi;\psi\})$ in \cite{cdf92} satisfying \eqref{I:phidphi} and \eqref{I:psidpsi} with $\wh{\phi}(0)=\wh{\tilde{\phi}}(0)=1$ and an associated
	biorthogonal wavelet filter bank $(\{\tilde{a};\tilde{b}\},\{a;b\})$ given by
	\begin{align*}
	 a=&\left\{\tfrac{1}{4},\tfrac{1}{2},\tfrac{1}{4}\right\}_{[-1,1]}, \quad b=\left\{-\tfrac{1}{8},-\tfrac{1}{4},\tfrac{3}{4},-\tfrac{1}{4},-\tfrac{1}{8}\right\}_{[-1,3]},\\
	\tilde{a}=&\left\{-\tfrac{1}{8}, \tfrac{1}{4}, \tfrac{3}{4}, \tfrac{1}{4}, -\tfrac{1}{8} \right\}_{[-2,2]}, \quad \tilde{b}=\left\{-\tfrac{1}{4}, \tfrac{1}{2}, -\tfrac{1}{4}\right\}_{[0,2]}.
	\end{align*}
	Note that $\phi$ is the hat function or the centered B-spline of order 2, $B_{2}$, which is defined as follows
	\be \label{def:hat}
	\phi(x)=B_{2}(x)=\begin{cases}
	1+x, &x\in [-1,0),\\
	1-x, &x\in [0,1],\\
	0, &\text{otherwise}.
	\end{cases}
	\ee
	By calculation, we have $\sm(\phi)=\sm(a)=1.5$, $\sm(\tilde{\phi})=\sm(\tilde{a}) \approx 0.440765$ and $\sr(a)=\sr(\tilde{a})=2$, where
$\sm(a)$ is defined in \cite[(5.6.44)]{hanbook}. Using Algorithm~\ref{alg:Phi} with $\pp(x)=x$, we have $n_{\phi}=3$ and
	\be \label{ex:legallphi}
	\phi^{L}=A_{c} \phi^{c} \quad \mbox{with} \quad
	\phi^{c}=\begin{bmatrix}
	\phi(\cdot-2)\\
	\phi(\cdot-1)\\
	\phi\chi_{[0,\infty)}
	\end{bmatrix},
	\quad A_{c}=\begin{bmatrix}
	1 & 0 & 0\\
	0 & 1 & 0\\
	\end{bmatrix}.
	\ee
	Setting $\pp(x)=x$ is what helps us to fulfill the homogeneous Dirichlet boundary condition. Similarly, using Algorithm~\ref{alg:Phi} with $\pp(x)=x$, we have $n_{\mathring{\phi}}=3$ and \eqref{ex:legallphi} with $\phi$ being replaced by $\mathring{\phi}$. Using Algorithm~\ref{alg:Phidual} with $m=\sr(\tilde{a})=2$ and $n_{\tilde{\phi}}=3$, we have
	\be \label{ex:legalltphi}
	\tilde{\phi}^{L}=\tilde{A}_{c} \tilde{\phi}^{c} \quad \mbox{with} \quad  \tilde{\phi}^{c}=
	\left[ \begin{matrix}
	\tilde{\phi}(\cdot-2)\\
	\tilde{\phi}(\cdot-1) \chi_{[0,\infty)}\\
	\tilde{\phi} \chi_{[0,\infty)}\\
	\tilde{\phi}(\cdot+1) \chi_{[0,\infty)}
	\end{matrix} \right], \quad
	\tilde{A}_{c} =
	\left[ \begin {array}{cccc} 1 & 0 & -1 & -2
	\\ \noalign{\medskip} 0 & 1 & 2 & 3
	\end {array} \right].
	\ee
    Setting $m=\sr(\tilde{a})=2$ is what helps us to achieve the maximum vanishing moments. Similarly, using Algorithm 3 with $m=\sr(\tilde{\mathring{a}})=2$ and $n_{\tilde{\mathring{\phi}}}=3$, we have \eqref{ex:legalltphi} with $\tilde{\phi}^{L}, \tilde{\phi}^{c}, \tilde{A}_{c}, \tilde{\phi}$ being replaced by $\tilde{\mathring{\phi}}^{L}, \tilde{\mathring{\phi}}^{c}, \tilde{\mathring{A}}_{c}, \tilde{\mathring{\phi}}$ respectively. Moreover, $\phi^{L}$, $\psi^{L}$, $\tilde{\phi}^{L}$, and $\tilde{\psi}^{L}$ satisfy
	\begin{align}\phi^{L} = & 2\begin{bmatrix}
	0 & 0\\
	\frac{1}{2} & \frac{1}{4}
	\end{bmatrix} \phi^{L}(2\cdot)
	+ 2\begin{bmatrix}
	\frac{1}{4}\\
	\frac{1}{4}
	\end{bmatrix} \phi(2\cdot-3)
	+ 2 \begin{bmatrix}
	\tfrac{1}{2}\\
	0
	\end{bmatrix} \phi(2\cdot-4)
	+ 2\begin{bmatrix}
	\frac{1}{4}\\
	0
	\end{bmatrix} \phi(2\cdot-5), \label{ex:legallphiL} \\
	\tilde{\phi}^{L}= & 2 \left[ \begin {array}{cc}
	0&-{\frac{1}{4}}\\ \noalign{\medskip}{
		\frac{1}{2}}&{\frac{3}{4}}\end {array} \right] \tilde{\phi}^{L}(2\cdot)
	+ 2 \left[ \begin {array}{c} {\frac{1}{4}}\\ \noalign{\medskip}\frac{1}{4}
	\end {array} \right] \tilde{\phi}(2\cdot-3)
	+ 2 \left[ \begin {array}{c} {\frac{3}{4}}\\ \noalign{\medskip}-{\frac{
			1}{8}}\end {array} \right] \tilde{\phi}(2\cdot-4)
	+ 2 \left[ \begin {array}{c} \frac{1}{4}\\ \noalign{\medskip}0\end {array}
	\right] \tilde{\phi}(2\cdot-5) \nonumber\\
	&+ 2 \left[ \begin {array}{c} -\frac{1}{8}\\ \noalign{\medskip}0\end {array}
	\right] \tilde{\phi}(2\cdot-6), \label{ex:legalltphiL}\\
	\psi^{L} = & 2\begin{bmatrix}
	b(0) & b(-1)\\
	-\frac{27}{64} & \frac{37}{128}
	\end{bmatrix} \phi^{L}(2\cdot)
	+ 2 \begin{bmatrix}
	b(1)\\
	\frac{1}{64}
	\end{bmatrix} \phi(2\cdot-3)
	+ 2\begin{bmatrix}
	b(2)\\
	\frac{5}{64}
	\end{bmatrix} \phi(2\cdot-4)
	+ 2 \begin{bmatrix}
	b(3)\\
	\frac{5}{128}
	\end{bmatrix} \phi(2\cdot-5), \nonumber\\
	\tilde{\psi}^{L}= & 2 \left[ \begin {array}{cc} -\frac{5}{32} & -\frac{3}{16}
	\\ \noalign{\medskip}-{\frac{1}{2}}& 1 \end {array}
	\right] \tilde{\phi}^{L}(2\cdot)
	+ 2 \begin{bmatrix}
	\frac{1}{2} \\
	0
	\end{bmatrix}\tilde{\phi}(2\cdot-3)
	+ 2 \begin{bmatrix}
	-\frac{1}{4} \\
	0
	\end{bmatrix}\tilde{\phi}(2\cdot-4). \label{ex:legalltpsiL}
	\end{align}
	Additionally, $\mathring{\phi}^{L}$, $\tilde{\mathring{\phi}}^{L}$ respectively satisfy \eqref{ex:legallphiL} and \eqref{ex:legalltphiL} with $\phi$, $\tilde{\phi}, \phi^{L}, \tilde{\phi}^{L}$ being properly replaced by $\mathring{\phi}$, $\tilde{\mathring{\phi}}, \mathring{\phi}^{L}$, $\tilde{\mathring{\phi}}^{L}$. On the other hand, $\mathring{\psi}^{L}$ and $\tilde{\mathring{\psi}}^{L}$ satisfy
	\begin{align}
	\mathring{\psi}^{L} = & 2\begin{bmatrix}
	\frac{2}{7} & 0\\
	-\frac{9}{28} & \frac{1}{4}
	\end{bmatrix} \mathring{\phi}^{L}(2\cdot)
	+ 2 \begin{bmatrix}
	-\frac{1}{2}\\
	-\frac{1}{16}
	\end{bmatrix} \mathring{\phi}(2\cdot-3)
	+ 2\begin{bmatrix}
	\frac{1}{7}\\
	\frac{5}{56}
	\end{bmatrix} \mathring{\phi}(2\cdot-4)
	+ 2 \begin{bmatrix}
	\frac{1}{14}\\
	\frac{5}{112}
	\end{bmatrix} \mathring{\phi}(2\cdot-5),  \label{ex:legallpsiLring}\\
	\tilde{\mathring{\psi}}^{L}= & 2 \left[ \begin {array}{cc} {\frac{5}{16}}&{\frac{3}{32}}\\
	-{\frac{1}{2}}&{\frac{5}{4}}\end {array}
	\right] \tilde{\mathring{\phi}}^{L}(2\cdot)
	+ 2 \begin{bmatrix} -{\frac{23}{32}}\\
	-{\frac{1}{4}} \end{bmatrix}
	\tilde{\mathring{\phi}}(2\cdot-3) + \begin{bmatrix} {\frac{23}{64}}\\
	{\frac{1}{8}} \end{bmatrix}
	\tilde{\mathring{\phi}}(2\cdot-4). \label{ex:legalltpsiLring}
	\end{align}
	Let $\cB_J=\Phi_J \cup \{\Psi_j \setsp j\ge J\}$ for all $J\ge 2$, where $\Phi_j$ and $\Psi_j$ are defined in \eqref{Phij} and \eqref{Psij} respectively, with $n_{\phi}=n_{\mathring{\phi}}=n_{\mathring{\psi}}=3$, $n_{\psi}=2$, $\phi^{R}=\mathring{\phi}^{L}(1-\cdot)$, and $\psi^{R}=\mathring{\psi}^{L}(1-\cdot)$. Let $\tilde{\cB}_J=\tilde{\Phi}_J \cup \{\tilde{\Psi}_j \setsp j\ge J\}$ for all $J \ge 3$, where $\tilde{\Phi}_j$ and $\tilde{\Psi}_j$ are defined in \eqref{tPhij} and \eqref{tPsij} respectively, with $n_{\tilde{\phi}}=n_{\tilde{\mathring{\phi}}}=n_{\tilde{\mathring{\psi}}}=3$, $n_{\tilde{\psi}}=2$, $\tilde{\phi}^{R}=\tilde{\mathring{\phi}}^{L}(1-\cdot)$, and $\tilde{\psi}^{R}=\tilde{\mathring{\psi}}^{L}(1-\cdot)$. According to \cite[Theorem 4.1]{hmw19} with $N=1$, $(\tilde{\cB}_J,\cB_J)$ form a biorthogonal Riesz basis of $L_{2}([0,1])$ for every $J \ge 3$. Note that $\vmo(\psi^L)=\vmo(\psi^R)=\vmo(\psi)=2=\sr(\tilde{a})$. However, $\vmo(\tilde{\psi}^{L})=\vmo(\tilde{\psi}^{R})=0$ and $\vmo(\tilde{\psi})=2=\sr(a)$. See Figure \ref{fig:legall} for the generators of $(\tilde{\cB}_{J},\cB_{J})$ for $J \ge 3$.
	
	As an implication of \cite[Theorem 4.1]{hmw19}, $(\tilde{\cB}_2,\cB_2)$ in Example~\ref{ex:legall} is also a biorthogonal Riesz basis of $L_{2}([0,1])$. Though $\cB_2$ is unchanged and is still given by \eqref{ex:legalltphiL} and \eqref{ex:legalltpsiL}, we point out that the boundary dual functions in $\tilde{\cB}_2$ may be different from \eqref{ex:legalltphiL}, \eqref{ex:legalltpsiL}, and \eqref{ex:legalltpsiLring}.
	
	It is indeed possible to replace the second component of $\psi^{L}$ with one that has a shorter support and symmetry. The same thing applies to $\psi^{R}=\mathring{\psi}^{L}(1-\cdot)$. However, our calculation suggests that such a choice leads to a much larger condition number, assuming the first components of $\psi^{L}$, $\mathring{\psi}^{L}$ are the same as our present choice.
\end{example}

\begin{figure}[htbp]
	\centering
	\begin{subfigure}[b]{0.24\textwidth}
		 \includegraphics[width=\textwidth]{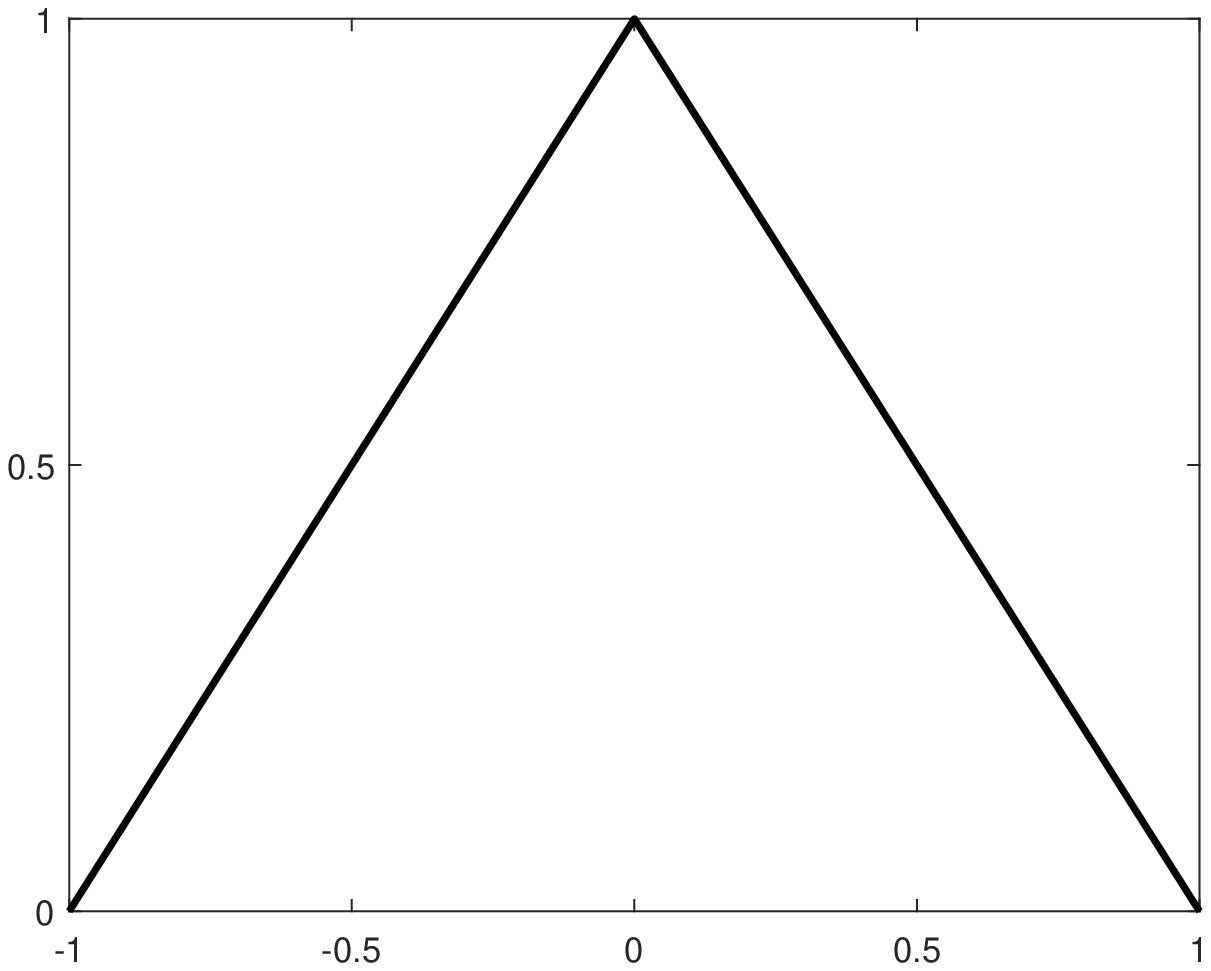}
		\caption{$\phi$}
	\end{subfigure}
	\begin{subfigure}[b]{0.24\textwidth}
		 \includegraphics[width=\textwidth]{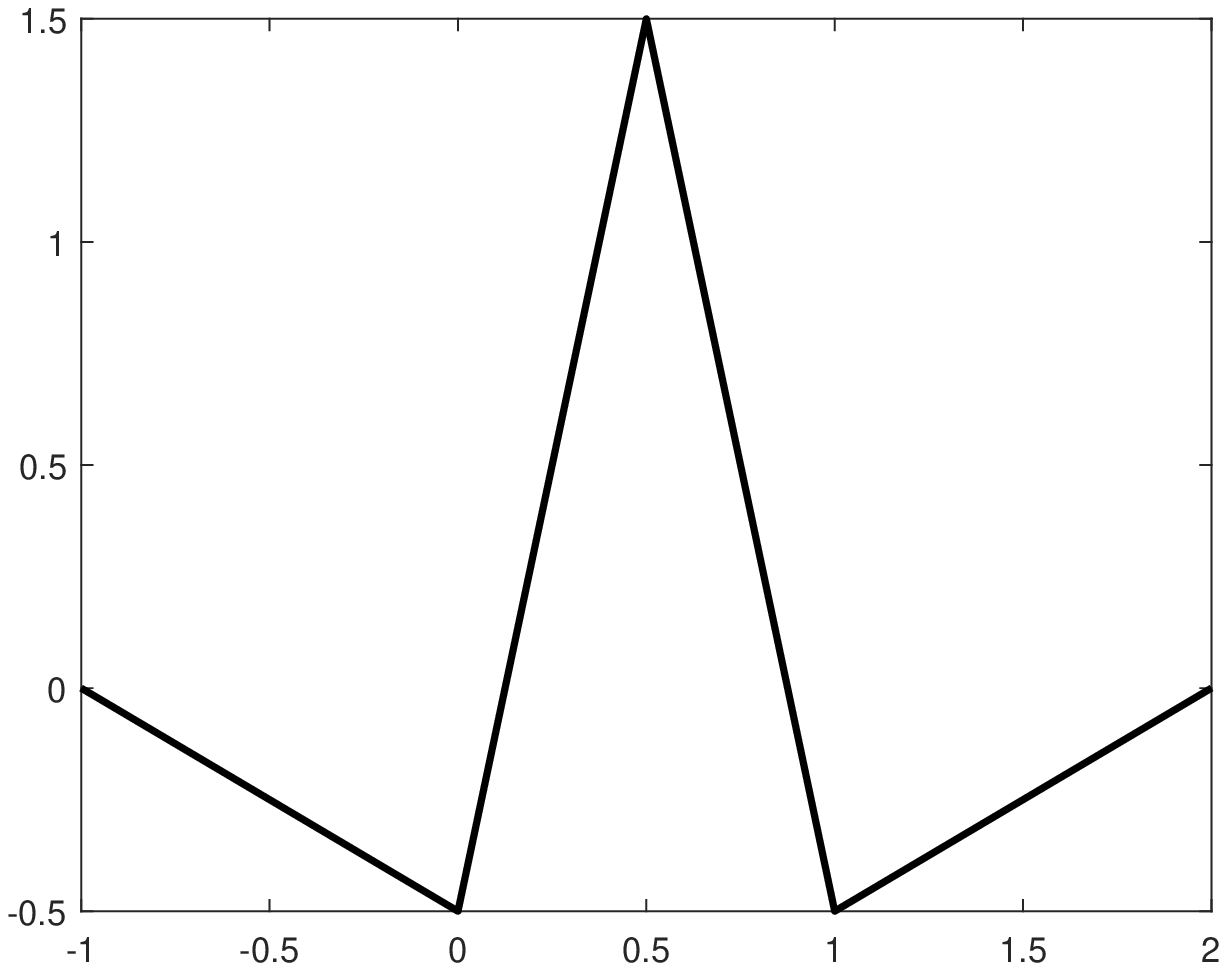}
		\caption{$\psi$}
	\end{subfigure}
	\begin{subfigure}[b]{0.24\textwidth}
		 \includegraphics[width=\textwidth]{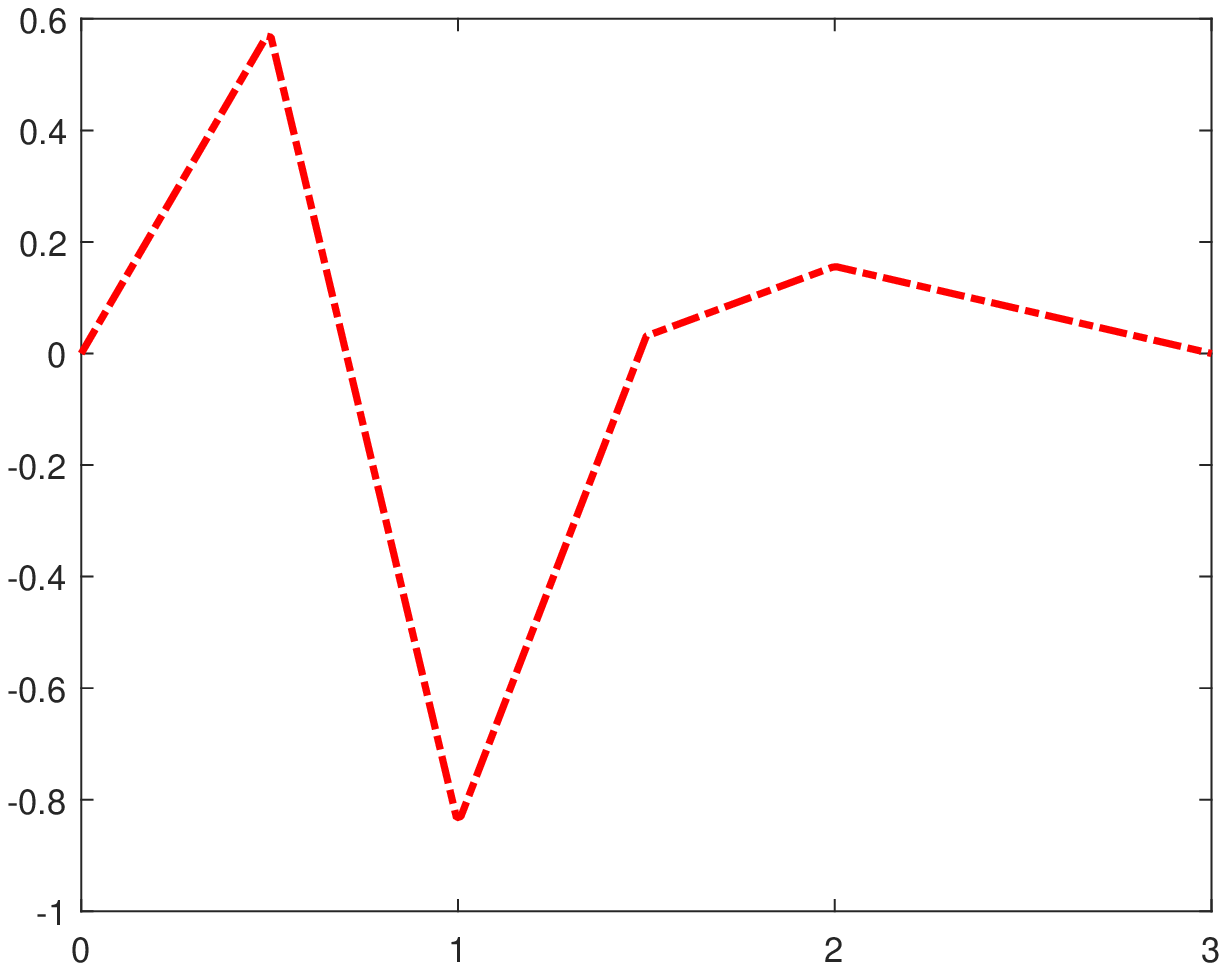}
		\caption{$\psi^{L}$}
	\end{subfigure}
	\begin{subfigure}[b]{0.24\textwidth}
		 \includegraphics[width=\textwidth]{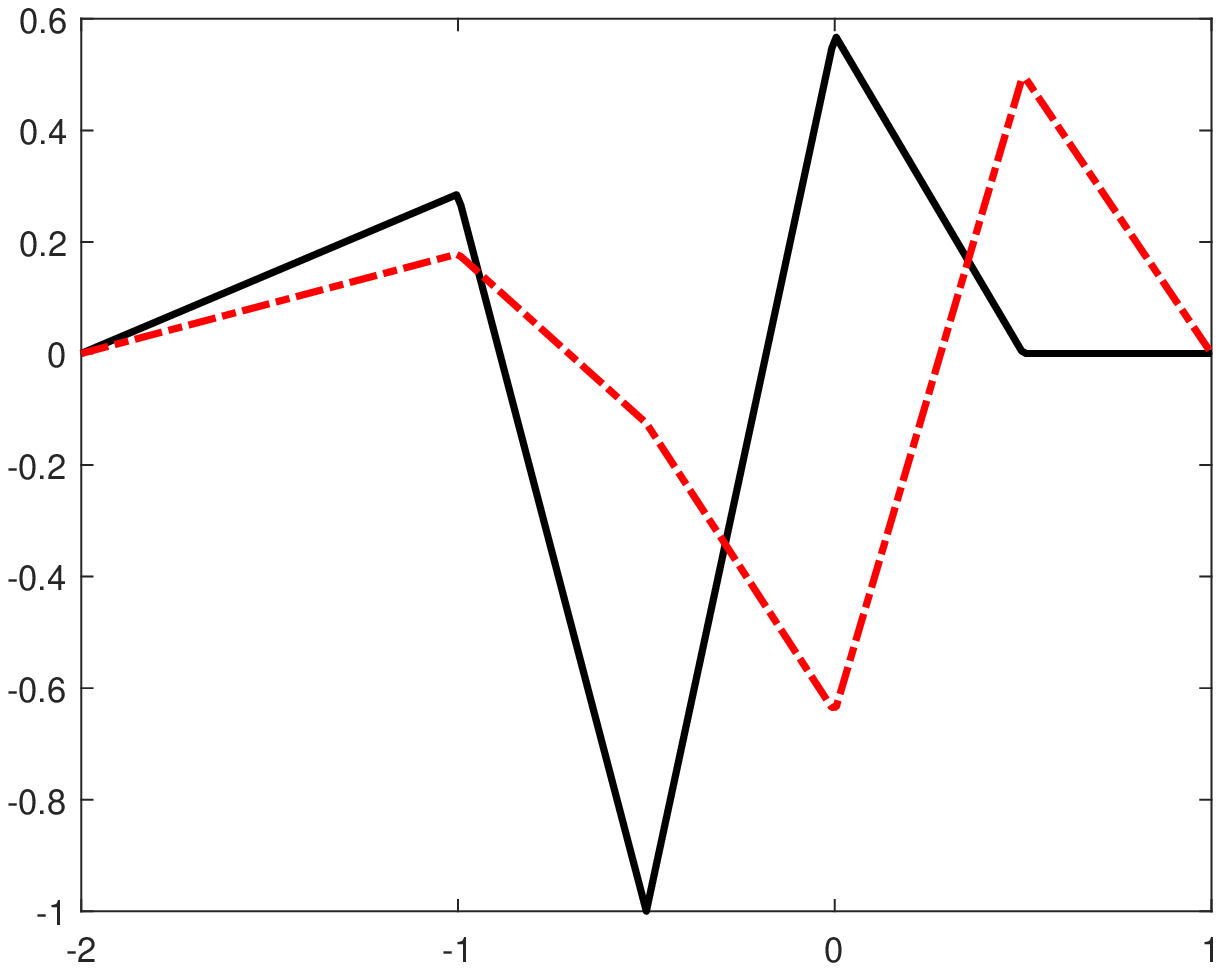}
		\caption{$\psi^{R}$}
	\end{subfigure}
	\begin{subfigure}[b]{0.24\textwidth}
		 \includegraphics[width=\textwidth]{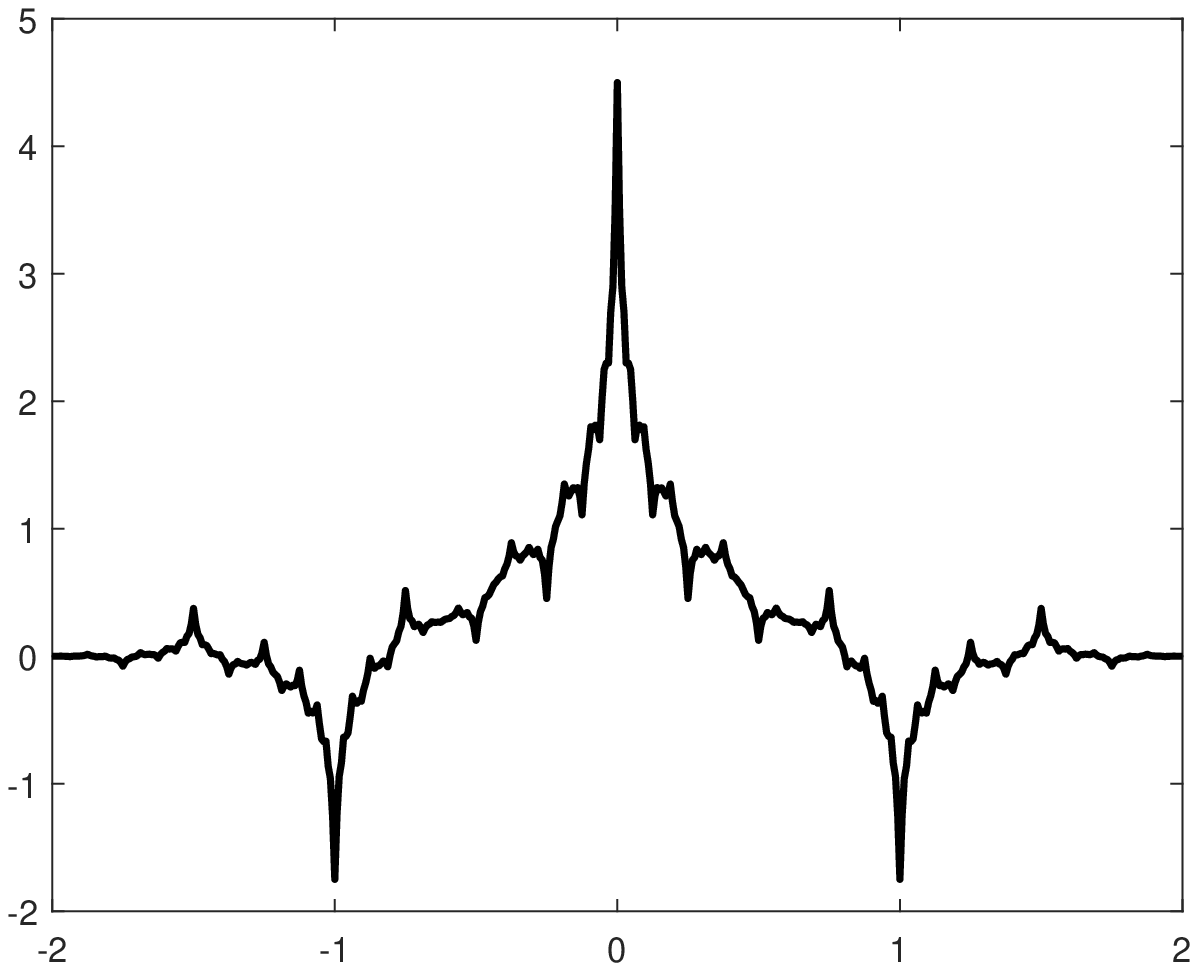}
		\caption{$\tilde{\phi}$}
	\end{subfigure}
	\begin{subfigure}[b]{0.24\textwidth}
		 \includegraphics[width=\textwidth]{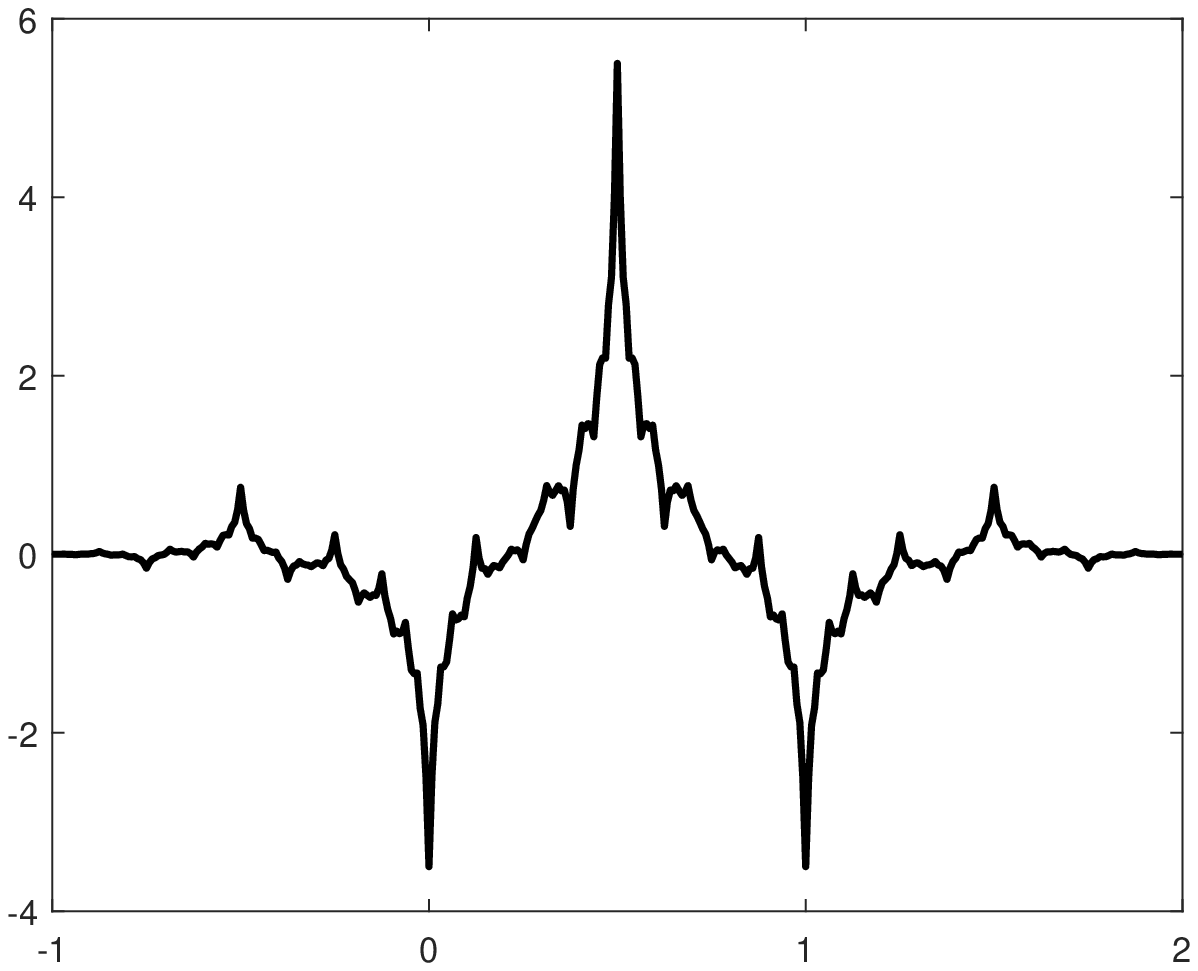}
		\caption{$\tilde{\psi}$}
	\end{subfigure}
	\begin{subfigure}[b]{0.24\textwidth}
		 \includegraphics[width=\textwidth]{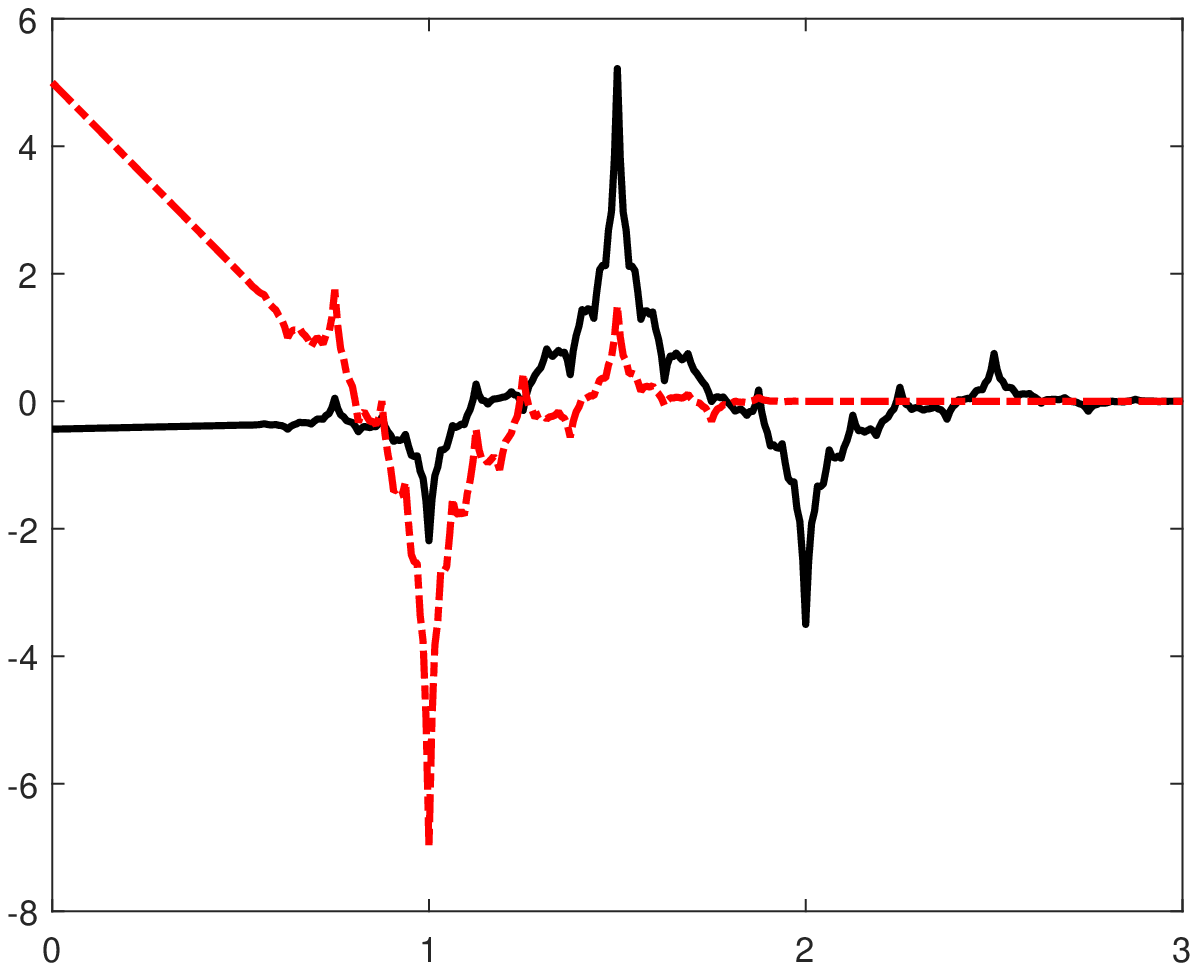}
		\caption{$\tilde{\psi}^{L}$}
	\end{subfigure}
	\begin{subfigure}[b]{0.24\textwidth}
		 \includegraphics[width=\textwidth]{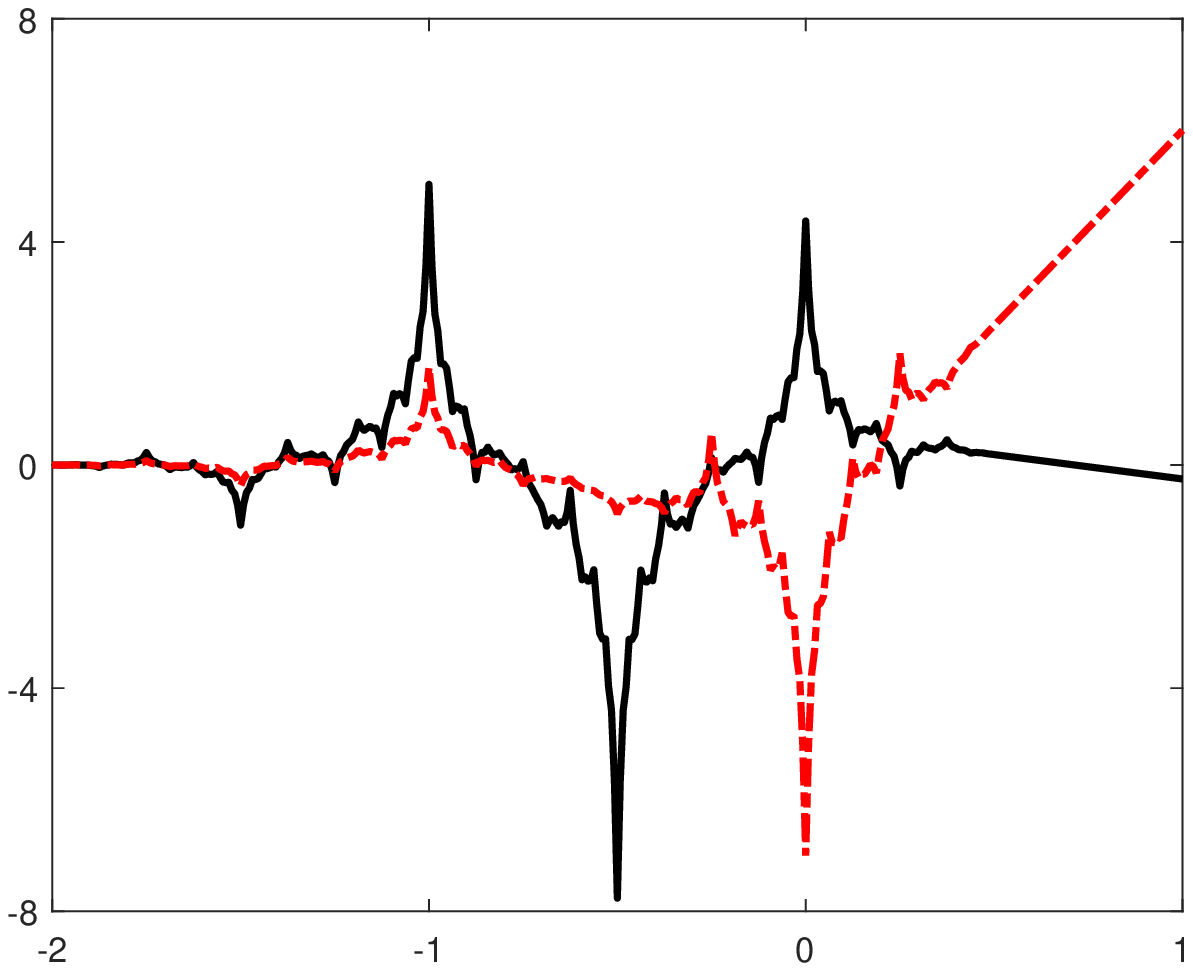}
		\caption{$\tilde{\psi}^{R}$}
	\end{subfigure}
	\caption{The generators of the biorthogonal wavelet $(\tilde{\cB}_{J},\cB_{J})$ of $L_{2}([0,1])$ with $J \ge 3$ in Example~\ref{ex:legall}. Black solid and red dashed lines correspond to the first and second components of the vector function respectively. (a) is the refinable function $\phi$. (b) is the wavelet function $\psi$. (c) is the second component of the left boundary wavelet function $\psi^{L}$ (the first component coincides with panel (b)). (d) is the right boundary wavelet function $\psi^{R}$. (e) is the dual refinable function $\tilde{\phi}$. (f) is the dual wavelet function $\tilde{\psi}$. (g) is the left boundary dual wavelet function $\tilde{\psi}^{L}$. (h) is the right boundary dual wavelet function $\tilde{\psi}^{R}$.}
	\label{fig:legall}
\end{figure}

\section{Numerical examples using Riesz wavelets in $L_2([0,1])$}

In this section, we consider two different model problems. One of which is the Helmholtz equation and the other is the biharmonic equation. Without loss of generality, we consider $\mathcal{I}=[0,1]$ as our domain of interest. 
The performance is measured by the relative $L_2$ error given by
\[
{\|e_{N}\|_{L_{2}}}/{\| u\|_{L_{2}}}:= {\|u_{N} -u\|_{L_{2}}}/{\| u\|_{L_{2}}} \qquad \mbox{with}\qquad e_N:=u_N-u,
\]
where $u$ is the true solution and $u_N$ is the numerically computed approximated solution with $N$ being its corresponding finest scale level for computing $u_N$. The meaning of $N$ will be further clarified below.
Because both $u$ and $u_N$ in our Helmholtz examples are highly oscillating functions, for the purpose of high accuracy and fair comparison, we symbolically compute both $\|u_N-u\|_{L_2}$ and $\|u\|_{L_2}$ with extremely high precision.

Let us first exemplify some advantages of wavelet-based methods over the standard finite element method by comparing the condition numbers of the following system $\cB_{2,N}$ truncated at the finest scale levels $N$ of the wavelet basis $\cB_2$ in $L_2([0,1])$ constructed in Example~\ref{ex:legall} of Section~3:
\be \label{B2N}
\cB_{2,N}:=\Phi_2 \cup \{\Psi_j \setsp 2\le j\le N-1\},\qquad N\ge 2
\ee
to the condition numbers of the standard finite element method using only the shifts of the refinable function $\phi=B_{2}$ defined in \eqref{def:hat} as follows:
\begin{equation}\label{hats}
\mbox{FEM}_N:=\{B_2(2^N\cdot-k) \setsp k=1,\dots,2^{N}-1\},\qquad N\ge 2.
\end{equation}
Note that these two systems $\cB_{2,N}$ and $\mbox{FEM}_N$
satisfy the homogeneous Dirichlet boundary conditions on $[0,1]$ and generate the same finite-dimensional space for all $N\ge 2$. See Table~\ref{table:cnlegall} for details.

\begin{center}
	\begin{tabular}{c c | c  c | c c}
		\hline
		Scale $N$ &Size
& Mass (FEM) & Stiffness (FEM) & Mass (Wavelet) & Stiffness (Wavelet) \\
		\hline
		11 &2047 & 3.0000 & $1.6999 \times 10^{6}$ & 18.4336 & 16.9644 \\
		12 &4095 & 3.0000 & $6.7929 \times 10^{6}$  & 19.2825 & 17.2715 \\
		13 &8191 & 3.0000 & $2.7198 \times 10^{7}$  & 20.0209 & 17.5118 \\
		14 &16383 & 3.0000 & $1.0879 \times 10^{8}$  & 20.6658 & 17.7025 \\
		\hline
	\end{tabular}
	\captionof{table}{Condition numbers of coefficient matrices using the finite element basis $\mbox{FEM}_N$ in \eqref{hats} and the wavelet basis $\cB_{2,N}$ in \eqref{B2N} of Example~\ref{ex:legall} with Dirichlet homogeneous boundary conditions. The mass and stiffness matrices have different normalizations. For the mass matrix, we normalize each element of  $\mbox{FEM}_N$ and $\cB_{2,N}$ such that it has norm equal to $1$. For the stiffness matrix, we normalize each element of $\mbox{FEM}_N$ and $\cB_{2,N}$ such that its derivative has norm equal to $1$.}
\label{table:cnlegall}
\end{center}


\subsection{Helmholtz equation}
Let $N\ge 2$ be the finest scale level for computation.
In this section, we shall modify the truncated system $\cB_{2,N}$ in \eqref{B2N} of the wavelet basis $\cB_2$ for $L_2([0,1])$ constructed in Example~\ref{ex:legall} of Section~3 to solve the Helmholtz equation.

Consider the following 1D model problem
\begin{align}
& -u''-k^2u=f \quad \text{on }\; (0,1), \label{model:helmholtz}\\
& u(0)=0, \quad u'(1)-ik u(1)=0, \label{model:helmholtzBC}
\end{align}
where $f \in L_{2}([0,1])$ and $\kappa >0$.
Since the generators $\phi,\tilde{\phi}$ for the biorthogonal wavelets $(\tilde{\cB}_2,\cB_2)$ in Example~\ref{ex:legall} has smoothness $\sm(\phi)=\sm(a)=1.5$ and $\sm(\tilde{\phi})=\sm(\tilde{a})\approx 0.440765$, according to Theorem~\ref{thm:rieszsobolev}, $\AS_2^\tau(\phi;\psi)$ must be a Riesz basis in the Sobolev space $H^\tau(\R)$ for all $-0.440765<\tau<1.5$.
Hence, after renormalizing each element in $\cB_{2,N}$ in $H^1$
and modifying the boundary wavelets to satisfy \eqref{model:helmholtzBC} (more details below),
then we can use the modified $\cB_{2,N}$ satisfying the boundary conditions in \eqref{model:helmholtzBC} to solve Helmholtz equation in \eqref{model:helmholtz}. The modification is admittedly heuristic by nature, but the modified $\cB_{2,N}$ still appears to be a Riesz basis in $H^{1}([0,1])$.

To capture the highly oscillating waves in the solution of the Helmholtz equation,
we additionally supplement our wavelet basis on $[0,1]$ with special waves. Suppose we have a non-overlapping partition of the unit interval; i.e., $[0,1] = \cup_{l=1}^{M} [a_{l},b_{l}]$, where $M$ denotes the number of partitions. Then, these special waves take the following form
\[
s^{+}_{l}(x):= (e^{ikx}-(\lambda^{+}_{1,l}+\lambda^{+}_{2,l}x))|_{[a_{l},b_{l}]}, \quad
s^{-}_{l}(x):= (e^{-ikx}-(\lambda^{-}_{1,l}+\lambda^{-}_{2,l}x))|_{[a_{l},b_{l}]},
\quad
l=1,\dots,M.
\]
Furthermore, the parameters $\lambda^{+}_{1,l}, \lambda^{+}_{2,l}, \lambda^{-}_{1,l}, \lambda^{-}_{2,l}$ are chosen so that
\begin{align*}
& s^{+}_{l}(a_{l})=s^{+}_{l}(b_{l})=s^{-}_{l}(a_{l})=s^{-}_{l}(b_{l})=0, \quad l=1,\dots,M-1, \\
& s^{+}_{M}(a_{M})=s^{-}_{M}(a_{M})=0, \quad (s^{+}_{M})'(1)-iks^{+}_{M}(1)=(s^{-}_{M})'(1)-iks^{-}_{M}(1)=0.
\end{align*}

Let us disregard the normalization constants for a moment, we note that our original right boundary elements are
\[
\mathcal{R}_{N}:= \{\phi(2^{2}\cdot-3)\} \cup \{\psi^{R}(2^{j}\cdot-(2^j-1)) \setsp 2 \le j \le N-1\}.
\]
Denote $\psi^{R,n}$ to be the $n$-th component of $\psi^{R}$. To ensure all right boundary elements satisfy the radiation boundary condition, we introduce a heuristic modification to a subset of these elements. Define
\[
\mathcal{S}_{N}:=\{\phi(2^{2}\cdot-3)\} \cup \{\psi^{R,2}(2^{j}\cdot-(2^j-1)) \setsp 2 \le j \le N-1\},
\]
where $\psi^{R}=\mathring{\psi}^{L}(1-\cdot)$ and $\mathring{\psi}^{L}$ is defined in \eqref{ex:legallpsiLring}. For each $g \in \mathcal{S}_{N}$, let $\text{supp}(g)$ be the support of $g$, $\text{supp}(g):=[l_{g},1]$, and define $g^{\text{new}}(x):=g(x)-(\lambda^{g}_{1}+\lambda^{g}_{2} x)$, where $\lambda^{g}_{1}, \lambda^{g}_{2} \in \C$. Instead of $\mathcal{R}_{N}$, we shall use the following modified right boundary elements
\[
\mathcal{R}_{N}^{\text{new}}:=\{\psi^{R,1}(2^{j}\cdot-(2^j-1)) \setsp j \ge 2\} \cup \{g^{\text{new}}\setsp g^{\text{new}}(l_{g})=0, (g^{\text{new}})'(1)-ikg^{\text{new}}(1)=0, g \in \mathcal{S}_N{}\}.
\]
To further reduce the condition number, we may first find a pair of linear combinations of the original right boundary wavelet functions at a given scale level, perform the modification, and find another pair of linear combinations of the modified right boundary wavelet functions at the same scale level. The same pairs of linear combinations are then applied to all scale levels.

Next, we normalize our modified wavelet basis on $[0,1]$ and special waves so that each of their first derivatives has norm equal to 1. Let $\{g_{n}\}_{n=1}^{2^{N}+2M-1}$ be the enumerated basis elements starting from the wavelet basis (ordered first by scales then by shifts) and ending with the special waves. The Galerkin formulation for \eqref{model:helmholtz} and \eqref{model:helmholtzBC} is
\[
\sum_{n=1}^{2^{N}+2M-1} \left( \langle g'_{l},g'_{n} \rangle-k^2\langle g_l,g_n\rangle-g_l(1)g_n'(1)\right) c_{n} = \langle g_{l},f \rangle, \quad l=1,\dots,2^{N}+2M-1.
\]

Consider the linear system of equations $Ax=b$ induced by the above Galerkin formulation with
\[
A=\begin{bmatrix}
A_{1} & A_{2}\\
A_{3} & A_{4}
\end{bmatrix}, \quad
x=\begin{bmatrix}
x_{1}\\
x_{2}
\end{bmatrix}, \quad
b=\begin{bmatrix}
b_{1}\\
b_{2}
\end{bmatrix}, \quad
\]
such that $A_{1} \in \C^{n_{1} \times n_{1}}$, $A_{2} \in \C^{n_{1} \times n_{2}}$, $A_{3} \in \C^{n_{2} \times n_{1}}$, $A_{4} \in \C^{n_{2} \times n_{2}}$, $x_{1},b_{1} \in \C^{n_{1}}$, and $x_{2},b_{2} \in \C^{n_{2}}$. Then, $(A_{1}-A_{2}A_{4}^{-1}A_{3}) x_{1} = b_{1} - A_{2}A_{4}^{-1}b_{2}$ and $A_{4} x_{2} = b_{2}-A_{3}x_{1}$. For the two following examples, we shall use two different conditions numbers $\kappa$ and $\kappa^{*}$. The former corresponds to the condition number of the entire coefficient matrix $A$; meanwhile, the latter corresponds to the condition number of $A_{1}-A_{2}A_{4}^{-1}A_{3}$. The submatrix $A_{4}$ coincides with the part of the coefficient matrix, where the special waves reside. The reason for having such a quantity is because generally the number of special waves is very small. Consequently, the matrix $A_{4}$ is often well-conditioned and can be inverted without any problem. Therefore, we can also treat the condition number of $(A_{1}-A_{2}A_{4}^{-1}A_{3})$ as a good indicator of how well-conditioned the entire coefficient matrix $A$ is and we may use it to reduce the size of $A$ with improved condition numbers.

To assess the performance of our enriched wavelet basis (the modified $\cB_{2,N}$ in \eqref{B2N} plus the special waves), we compare it with the standard wavelet-based Galerkin method (only the modified $\cB_{2,N}$ in \eqref{B2N}) and the (pollution free) finite difference method proposed in \cite{ww14}. Additionally, we utilize two distinct error measurements for the finite difference method: discrete and interpolation errors. Let $h$ be the grid size and $H_{\text{FD}}=1/h$. Define the discrete (relative) error as
\[
\text{Dis. Err.}:= \left. \left(\sum_{n=1}^{H_{\text{FD}}} |u(n h)-U_{H_{\text{FD}}}(n)|^{2}\right)^{\frac{1}{2}} \middle/ \left(\sum_{n=1}^{H_{\text{FD}}} |u(n h)|^{2}\right)^{\frac{1}{2}}, \right.
\]
where $u$ is the exact solution and $U_{H_{\text{FD}}}$ is the solution obtained from the finite difference method. The reason for introducing another measure of errors is because in a finite difference scheme, we do not know the values in between the discrete approximated solution. To partially reconcile the $L_{2}$ relative error of the wavelet-based Galerkin method and the discrete error of the finite difference method, we simply interpolate the discrete points obtained in the finite difference method. To calculate the interpolation error, we first interpolate the discrete points $U_{H_{\text{FD}}}$ with the hat function
\[
u_{N_{\text{FD}}}^{\text{int}}:=\sum_{n=1}^{H_{\text{FD}}} U_{H_{\text{FD}}}(n) B_{2}(2^{N_{\text{FD}}}\cdot-n) \qquad \text{with} \qquad N_{\text{FD}}:=\log_{2}(H_{\text{FD}}),
\]
and then calculate the relative $L_{2}$ error
\[
\|e_{N_{\text{FD}}}^{\text{int}}\|_{L_{2}}/{\| u\|_{L_{2}}}:= {\|u_{N_{\text{FD}}}^{\text{int}} -u\|_{L_{2}}}/{\| u\|_{L_{2}}} \qquad \mbox{with}\qquad e_{N_{\text{FD}}}^{\text{int}}:=u_{N_{\text{FD}}}^{\text{int}}-u.
\]

\begin{figure}[htbp]
	\centering
	\begin{subfigure}[b]{0.32\textwidth}
		 \includegraphics[width=\textwidth]{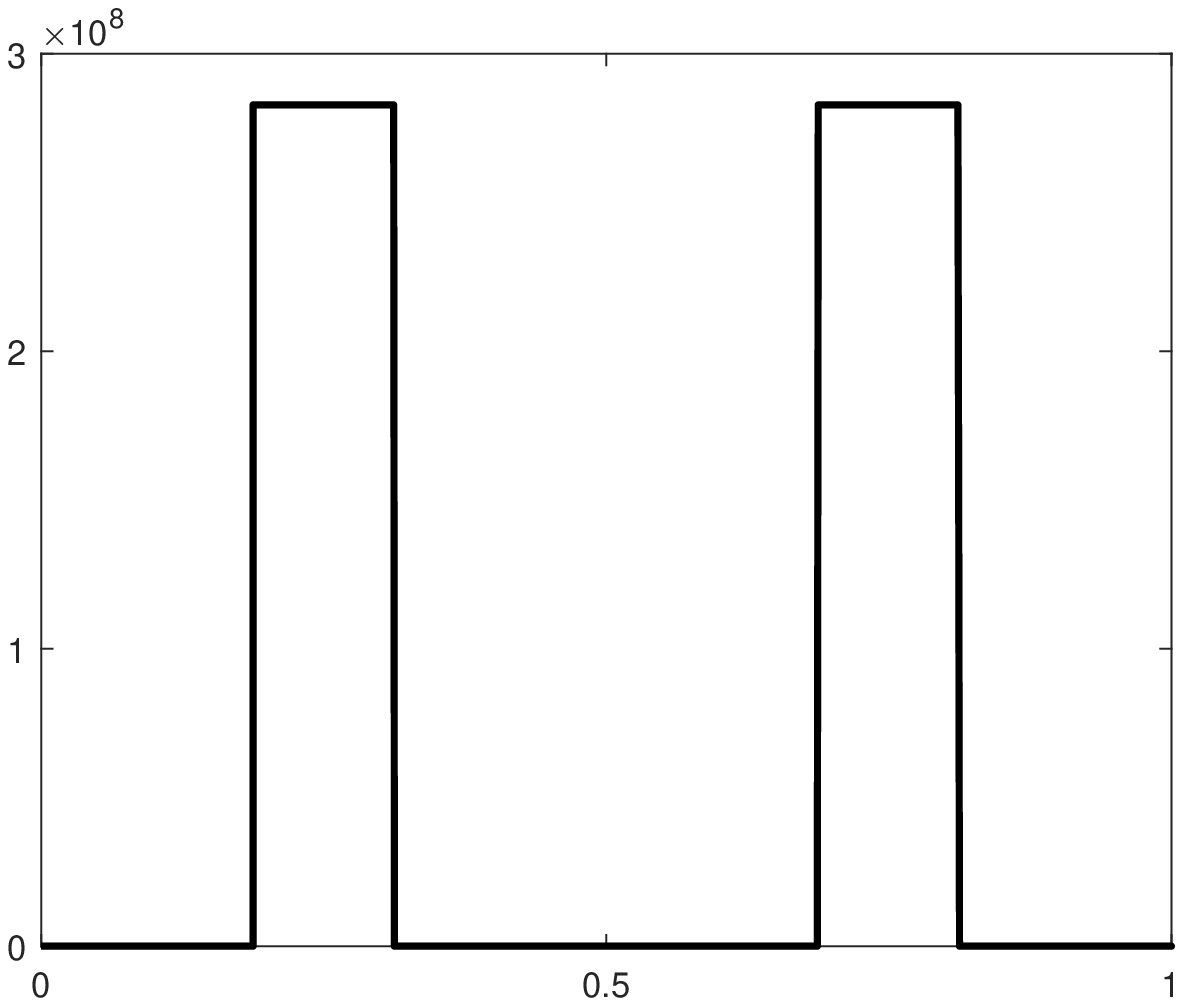}
	\end{subfigure}
	\begin{subfigure}[b]{0.32\textwidth}
		 \includegraphics[width=\textwidth]{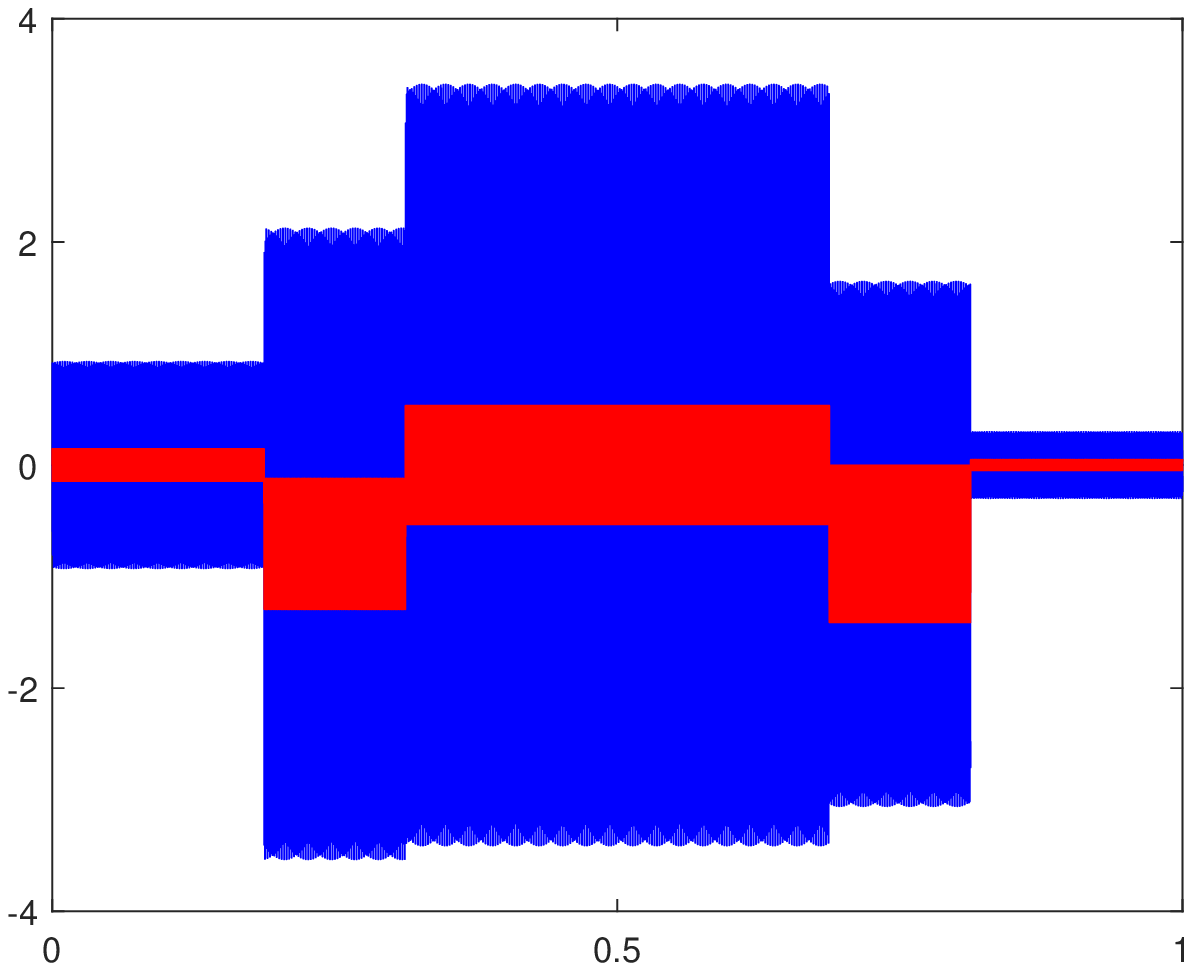}
	\end{subfigure}
	\begin{subfigure}[b]{0.32\textwidth}
		 \includegraphics[width=\textwidth]{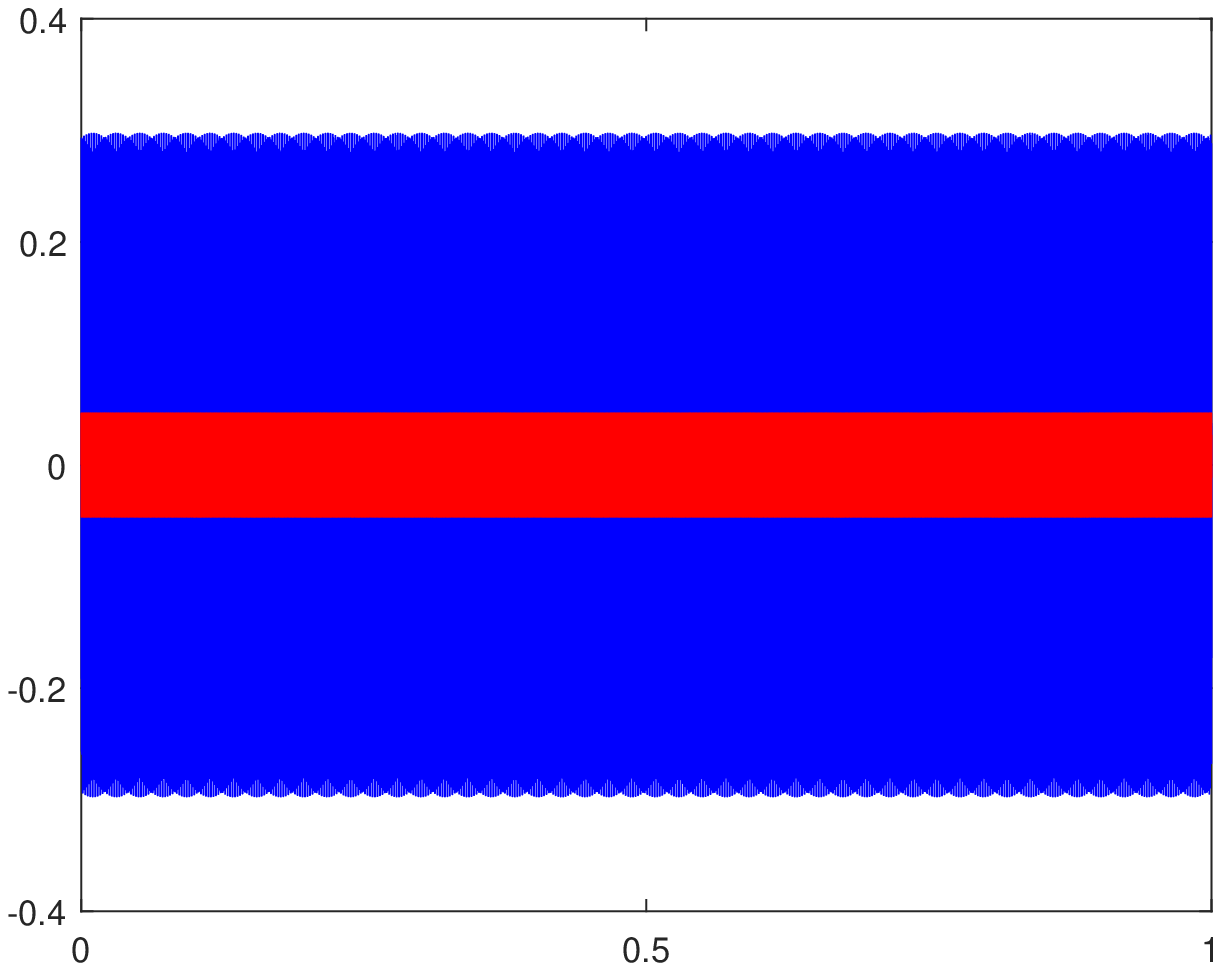}
	\end{subfigure}
	\caption{Plots for Example~\ref{ex:indicator}. The solid block-like plots are due to the high wave number, $k=20000$. Left: source term $f$. Middle: Real solutions obtained from finite difference method with a grid size $2^{-13}$ (blue) and wavelet-based Galerkin method with 10 special waves at scale level $N=4$ (red). Right: Imaginary solutions obtained from finite difference method with a grid size $2^{-13}$ (blue) and wavelet-based Galerkin method with 10 special waves at scale level $N=4$ (red).}
	\label{fig:ex:indicator}
\end{figure}

\begin{example} \label{ex:indicator}
	\normalfont
	Suppose $f:=2\sqrt{2} \times 10^{8} \times (\chi_{[\frac{3}{16},\frac{5}{16}]} + \chi_{[\frac{11}{16},\frac{13}{16}]})$ and $k=20000$ in \eqref{model:helmholtz}. The unit interval is partitioned into $[0,\frac{3}{16}] \cup [\frac{3}{16},\frac{5}{16}] \cup [\frac{5}{16},\frac{11}{16}] \cup [\frac{11}{16},\frac{13}{16}] \cup [\frac{13}{16},1]$. Table~\ref{table:ex:indicator} lists sizes of coefficient matrices, condition numbers, and errors. We see that the relative error corresponding to the wavelet-based Galerkin method with special waves hits $0\%$ at the scale level $N=4$, which implies that the solution is recovered exactly. See Figure~\ref{fig:ex:indicator} for plots of the source term as well as real and imaginary approximated solutions.
\end{example}

\begin{table}[htbp]
\begin{tabular}{c c|c c|c c||c c c c}
		\hline
		&& \multicolumn{2}{|c|}{Wavelet + 10 special waves} &  \multicolumn{2}{|c||}{Wavelet only} & \multicolumn{4}{|c}{Finite difference method in \cite{ww14}}\\
		\hline
		Scale $N$ & Size & $\kappa$ ($\kappa^{*}$) & $\frac{\|e_{N}\|_{L_{2}}}{\| u\|_{L_{2}}}$ & $\kappa$ & $\frac{\|e_{N}\|_{L_{2}}}{\| u\|_{L_{2}}}$ & Size & $\kappa$
&Dis. Err. & $\frac{\|e_{N_{\text{FD}}}^{\text{int}}\|_{L_{2}}}{\| u\|_{L_{2}}}$\\
		\hline
		3&	7& 124515 (9749) & 353.73\% & 4.673 & 73.54\% & 512& 1024 & 119.92\% & 106.00\%\\
		4&	15& 59241 (14372) & 0\% & 7.198 & 73.26\% & 1024& 3792 &   52.40\% & 85.97\%\\
		5&  31& 61083 (17585) & 0\% & 9.672 & 70.95\% & 2048 & 15336 &  861.35\% & 538.08\% \\
		6&	63& 62120 (18268) & 0\% & 11.716 & 69.59\% & 4096 & 8090 &  124.93\% & 159.72\%\\
		7&	127& 62679 (18534) & 0\% & 13.484 & 68.90\% & 8192& 29519 &  405.90\% & 273.96\%\\
		\hline
	\end{tabular}
	\caption{Error summary for Example~\ref{ex:indicator}. Size in the far left column does not include the 10 special waves. For the wavelet-based Galerkin method, $\kappa$ and $\kappa^{*}$ after we apply a diagonal preconditioner to the coefficient matrix so that each diagonal entry has modulus equal to 1.}
	\label{table:ex:indicator}
\end{table}

\begin{figure}[htbp]
	\centering
	\begin{subfigure}[b]{0.32\textwidth}
		 \includegraphics[width=\textwidth]{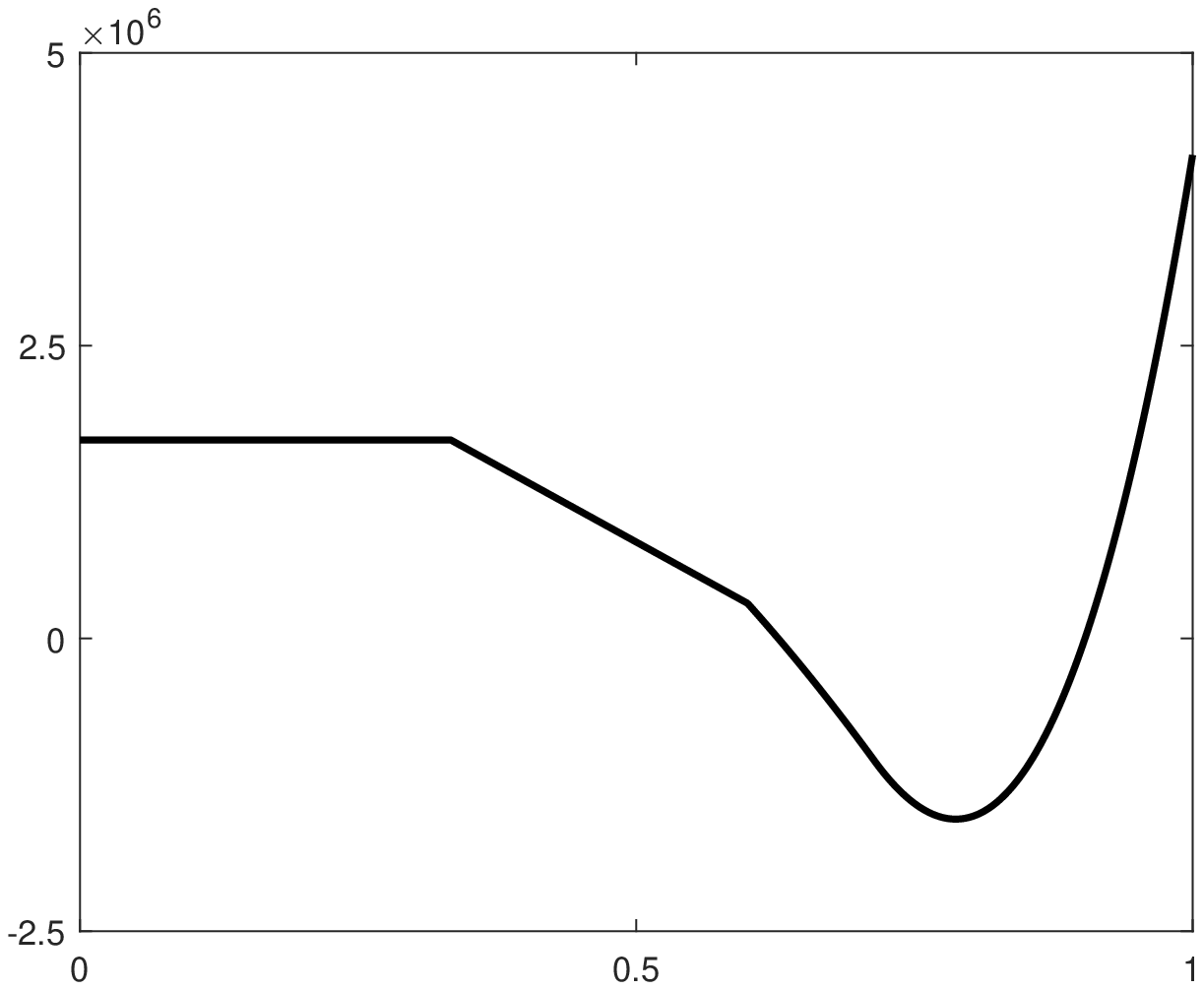}
	\end{subfigure}
	\begin{subfigure}[b]{0.32\textwidth}
		 \includegraphics[width=\textwidth]{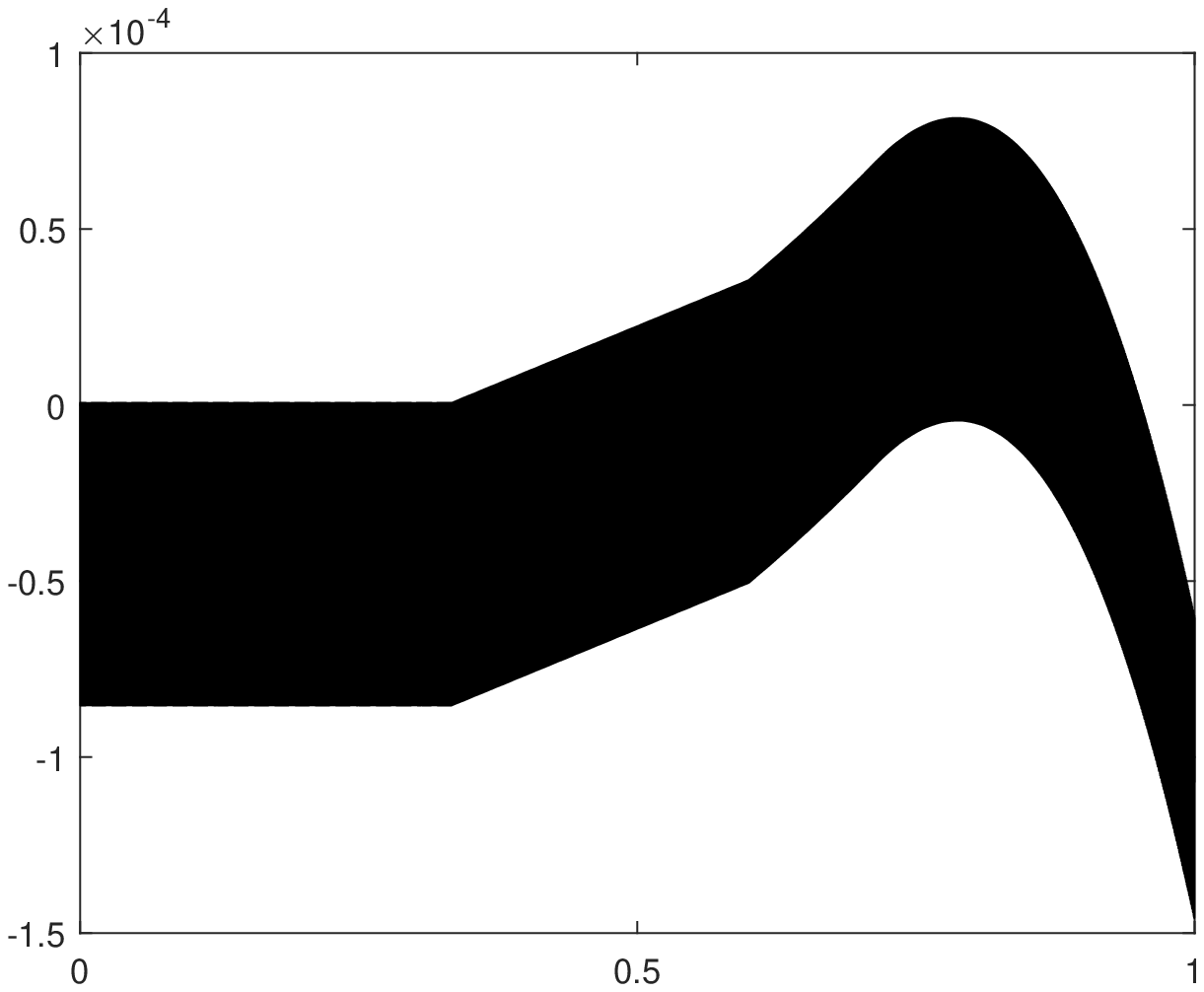}
	\end{subfigure}
	\begin{subfigure}[b]{0.32\textwidth}
		 \includegraphics[width=\textwidth]{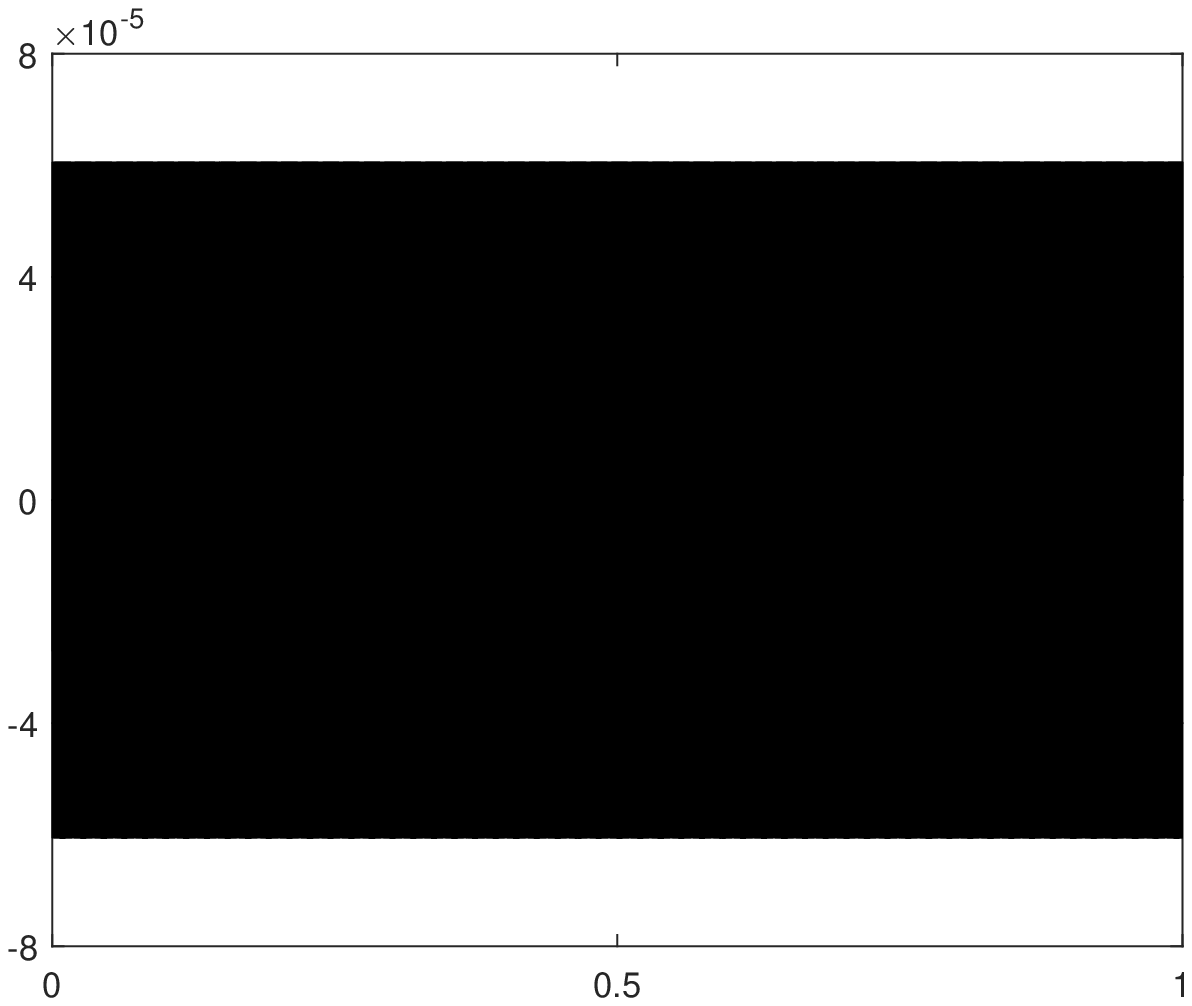}
	\end{subfigure}
	\caption{Plots for Example~\ref{ex:fpieceC1}. The solid block-like plots are due to the high wave number, $k=200000$. Left: source term $f$. Middle: Real part of the true solution. Right: Imaginary part of the true solution.}
	\label{fig:ex:fpieceC1}
\end{figure}

\begin{example} \label{ex:fpieceC1}
	\normalfont
	Suppose $f:=3000(\frac{2824}{5} \chi_{[0,\frac{1}{3}]} + (-1743+\frac{5729}{5}x) \chi_{[\frac{1}{3},\frac{3}{5}]} + (-3525x^{2}+700x+949) \chi_{[\frac{3}{5},\frac{5}{7}]} + (37000x^{3} - \frac{9370211}{175}x^2+\frac{3803172}{245}x+2397)\chi_{[\frac{5}{7},1]})$ and $k=200000$ in \eqref{model:helmholtz}. The unit interval is partitioned into $[0,\frac{3}{8}] \cup [\frac{3}{8},\frac{5}{8}] \cup [\frac{5}{8},1]$. Table~\ref{table:ex:fpieceC1} lists sizes of coefficient matrices, condition numbers, and errors. See Figure~\ref{fig:ex:fpieceC1} for plots of the source term as well as the real and imaginary true solutions.
\end{example}

\begin{table}
	\begin{tabular}{c c|c c|c c||c c c c}
		\hline
		&& \multicolumn{2}{|c|}{Wavelet + 6 special waves} &  \multicolumn{2}{|c||}{Wavelet only} & \multicolumn{4}{|c}{Finite difference method in \cite{ww14}}\\
		\hline
		Scale $N$ & Size & $\kappa$ ($\kappa^{*}$) & $\frac{\|e_{N}\|_{L_{2}}}{\| u\|_{L_{2}}}$ & $\kappa$ & $\frac{\|e_{N}\|_{L_{2}}}{\| u\|_{L_{2}}}$ & Size & $\kappa$
&Dis. Err. & $\frac{\|e_{N_{\text{FD}}}^{\text{int}}\|_{L_{2}}}{\| u\|_{L_{2}}}$\\
		\hline
		6&	63& 282841 (72458) &  0.233\%&	 11.715 & 83.50\% & 512 &	 1223&	3.749\% & 111.61\%\\
		7&  127& 285432 (73951) &  0.063\%&	 13.484 & 83.09\% & 1024&	 4865&   3.244\% & 115.40\%\\
		8& 255& 286785 (75179) &  0.017\%&	 15.008 & 82.88\% & 2048&  19627&   2.978\% & 95.87\% \\
		9& 511&	287541 (76024) &  0.005\%&	 16.321 & 82.77\% & 4096&   8010&   0.415\% & 108.09\%\\
		10& 1023& 288105 (76805) &  0.002\%& 17.455& 82.72\% & 8192&  28748&  0.317\% & 114.42\%\\
		\hline
	\end{tabular}
	\caption{Error summary for Example~\ref{ex:fpieceC1}. Size in the far left column does not include the 6 special waves. For the wavelet-based Galerkin method, $\kappa$ and $\kappa^{*}$ after we apply a diagonal preconditioner to the coefficient matrix so that each diagonal entry has modulus equal to 1.}
	\label{table:ex:fpieceC1}
\end{table}

The condition numbers of coefficient matrices coming from the enriched wavelet basis with special waves are large. Our wavelet basis is definitely not the source of this problem, as its condition numbers are small (see $\kappa$ in Tables~\ref{table:ex:indicator} and \ref{table:ex:fpieceC1} under `Wavelet only'), but rather the large condition numbers are due to the enrichment by special waves. This phenomenon has been reported in the literature of numerical methods for the Helmholtz equation, particularly those that use the homogeneous solution of the underlying DE in the trial and test spaces \cite{hmp16}.

\subsection{Biharmonic equation}
In this subsection, we shall apply Example~\ref{ex:hmtcubdow} of Section~2 to solve the following biharmonic equation
\begin{align}
& u^{(4)}=f \quad \text{on} \quad (0,1) \label{biharmonic}\\
& u(0)=u(1)=u'(0)=u'(1)=0, \label{biharmonicBC}
\end{align}
where $f \in L_{2}([0,1])$. To justify our construction on $[0,1]$, we observe $\{\mathring{\phi};\mathring{\psi}\}$ is an orthogonal wavelet in $L_{2}(\R)$, where $\mathring{\psi}:=\text{diag}(\sqrt{3}/24,1/8) \psi''$ and $\mathring{\phi}:=(\chi_{[0,1]}, (\sqrt{3}(2x-1)\chi_{[0,1]})^\tp$ such that $\{\mathring{\phi}(\cdot-k)\}_{k\in \Z}$ and
$\{\phi''(\cdot-k)\}_{k\in \Z}$ generate the same shift-invariant space. Using the construction in Section~3, we obtain an orthogonal wavelet in $L_{2}([0,1])$ derived from $\{\mathring{\phi};\mathring{\psi}\}$. As a consequence, we can deduce a Riesz wavelet in $H^{2}([0,1])$ derived from $\{\phi;\psi\}$. Unlike in Section~3.1 with the Helmholtz equation, we choose our coarsest scale level to be $J=1$. That is, define
\begin{align}
&\Phi_{1}:=\{2^{-1/2}\phi_{1}(2\cdot-1),2^{-1/2}\phi_{2}(2\cdot-1)\}, \label{Phibihar}\\ &\Psi_{j}:=\{2^{-j/2}\psi_{1}(2^{j}\cdot-k),2^{-j/2}\psi_{2}(2^{j}\cdot-k) \setsp 0\le k \le 2^{j}-1\}, \quad j \ge 1, \label{Psibihar}\\
&\cB_{1}:= \Phi_{1} \cup \{\Psi_{j} \setsp j \ge 1\}. \nonumber
\end{align}
After an appropriate renormalization in $H^{2}$, the system $\cB_{1}$ forms a Riesz basis in $H^{2}([0,1])$. As before, assuming that $N$ is our finest scale level, we have the truncated system:
\[
\cB_{1,N}:= \Phi_{1} \cup \{\Psi_{j}: 1 \le j \le N-1\}, \quad N \ge 1,
\]
where $\Phi_{1}, \Psi_{j}$ are defined in \eqref{Phibihar} and \eqref{Psibihar}, respectively. Let $\{g_{n}\}_{n=1}^{2^{N+2}-2}$ be the ordered normalized wavelet basis $\cB_{1,N}$ (first by scales then by shifts). See \cite[Section 4]{hm19} for the explicit ordering and normalization. The Galerkin formulation for \eqref{biharmonic} and \eqref{biharmonicBC} is
\[
\sum_{n=1}^{2^{N+2}-2} \langle g_{l}'',g_{n}'' \rangle c_{n} = \langle g_{l},f\rangle, \quad l=1,\dots,2^{N+2}-2.
\]

\begin{example} \label{biharmonic:sin}
	\normalfont
	Suppose $f:=-6250000\pi^{2}((-\frac{3}{625}+(x^2-x)\pi^2)\sin(50\pi x) -\frac{4}{25}\pi\cos(50\pi x)(x-\frac{1}{2}))$ in \eqref{biharmonic}. The true solution is $u(x)=\sin(50\pi x)(-x^2+x)$. Table~\ref{table:ex:biharmonic} lists sizes of coefficient matrices, condition numbers, and errors. As can be seen from the table, the condition number of the coefficient matrix is identically equal to 1. I.e., the coefficient matrix is an identity matrix. Therefore, we do not need to solve any linear systems.
\end{example}

\begin{table}
	\begin{tabular}{c c c c c}
		\hline
		Scale $N$ & Size & $\kappa$ & $\frac{\|e_{N}\|_{L_{2}}}{\| u\|_{L_{2}}}$ & $\log_{2} \frac{\|e_{N-1}\|_{L_{2}}}{\|e_{N}\|_{L_{2}}}$ \\
		\hline
		6&	254 &	  1&  $3.803 \times 10^{-1}\%$ & -- \\
		7&	510 &	  1&  $2.369 \times 10^{-2}\%$ & 4.005\\
		8&	1022 &	  1&  $1.479 \times 10^{-3}\%$ & 4.001\\
		9&	2046 &	  1&  $9.244 \times 10^{-5}\%$ & 4.000\\
		10&	4094 &    1&  $5.778 \times 10^{-6}\%$ & 4.000\\
		\hline
	\end{tabular}
	\caption{Error summary for Example~\ref{biharmonic:sin}. The condition numbers of coefficient matrices are identically equal to $1$. $\log_{2} \frac{\|e_{N-1}\|_{L_{2}}}{\|e_{N}\|_{L_{2}}}$ is the convergence rate, which coincides with the sum rule order $\sr(a)$ in Example~\ref{ex:hmtcubdow} of Section~2.}
	\label{table:ex:biharmonic}
\end{table}

We conclude this section with one remark. Notice that we have deliberately used a derivative-orthogonal Riesz wavelet only for the biharmonic equation. There certainly exists a first-order derivative-orthogonal Riesz wavelet generated from the hat function, whose stiffness matrix has a condition number identically equal to 1. See \cite[Example 3.1]{hm19}. However, such a wavelet is a Riesz basis only in $H^{1}(\R)$, but not in $L_{2}(\R)$. On the contrary, the wavelet studied in Example~\ref{ex:legall} of Section~3 is a Riesz basis in both $L_{2}(\R)$ and $H^{1}(\R)$. To ensure that the condition numbers of mass and stiffness matrices are uniformly bounded, the wavelet studied in Example~\ref{ex:legall} of Section~3, which is a Riesz basis in both $H^1(\R)$ and $L_2(\R)$ after renormalization,
is a more suitable choice for solving the Helmholtz equation with very large wave numbers.

\section{Conclusion}
In summary we have reviewed some important research directions in the context of wavelet-based methods for numerical DEs. Two main results revisited were the construction of $m$-th order derivative orthogonal Riesz wavelets in the Sobolev space $H^{m}(\R)$ as well as a general construction of wavelets in $L_{2}([0,N])$ satisfying maximum vanishing moments and given boundary conditions. Furthermore, we pointed out that adaptive wavelet algorithms have received a lot of interest over the years, even though we did not delve into this topic in this paper. Earlier on, we also discussed the advantages that wavelets bring in solving numerical DEs and demonstrated their ability in solving some model problems. One thing worth emphasizing is that we do not claim wavelets can efficiently solve all problems in numerical DEs. However, for selected problems, wavelets may serve as an effective approach.

Going forward, there are a few interesting problems to consider. Developing wavelets with simple structures that are capable of effectively tackling high dimensional problems remains as a critical research direction. We restrict the domain of interest to $[0,1]^d$, since the construction of wavelets on a general bounded domain poses a great challenge. Currently, the existence of derivative-orthogonal Riesz wavelet in $\R^{d}$, where $d \ge 2$, is unknown. If such wavelets existed, then we would encounter a desired situation like Example~\ref{biharmonic:sin} in a multidimensional setting. I.e., we do not need to solve any linear systems for basic Poisson and biharmonic problems. Otherwise, we need to explore the possibility of finding Riesz wavelets in $\R^{d}$ with near orthogonality in their derivatives. The Helmholtz equation serves as a good motivation to study the construction of wavelets on a bounded interval satisfying mixed (Robin) boundary conditions and maximum vanishing moments. Unfortunately, we do not have any theoretical justification for the boundary elements modification we adopted in Section~4.1. One difficulty that arises from mixed boundary conditions is that we now have to deal with nonstationary wavelet systems, since the boundary elements are no longer scale invariant. We are curious to see how we can extend our construction in Section~3 to accommodate mixed boundary conditions. Still on the topic of the Helmholtz equation, we earlier observe that wavelets combined with special waves yield encouraging results. However, the corresponding condition numbers are still large and they will continue to grow in high dimensions. Naturally, we want to know if we can design a better wavelet-based method that captures both low and high frequencies more efficiently.

\end{document}